\documentclass[12pt, english, a4paper]{amsart}
\usepackage[margin=2.5cm]{geometry}
\usepackage{amssymb,amsmath,amsthm}
\usepackage{mathtools}
\usepackage{mathrsfs}
\usepackage{accents}
\usepackage[utf8]{inputenc}
\usepackage[T1]{fontenc}

\usepackage{caption}
\usepackage{subcaption}

\usepackage{xfrac}
\usepackage{enumerate}
\usepackage{graphics, float}

\usepackage[colorlinks]{hyperref}
\hypersetup{allcolors = black}

\usepackage{url}

\usepackage{svg}

\newtheorem{theorem}{Theorem}

\newtheorem{corollary}{Corollary}
\newtheorem{proposition}{Proposition}
\newtheorem{lemma}{Lemma}
\newtheorem{problem}{Problem}

\theoremstyle{definition} 
\newtheorem{definition}{Definition}

\theoremstyle{remark}
\newtheorem{remark}{Remark}

\DeclareMathOperator{\conv}{\mathrm{conv}\,}

\newcommand{\R}{\mathbb{R}}
\newcommand{\K}{\mathcal{K}}

\newcommand{\la}{\langle}
\newcommand{\ra}{\rangle}

\newcommand{\bd}{\partial}
\newcommand{\eps}{\varepsilon}

\newcommand{\dd}{\mathrm{d}}

\newcommand{\per}{P}

\newcommand\etfvar{w}
\newcommand\etfind{N}

\title{Large signed sums and the polarization constant of convex bodies}
\author{Gergely Ambrus}
\thanks{Research of G. A. was partially supported by the ERC Advanced Grant "GeoScape" no.  882971, by Hungarian National Research (NKFIH) grants no. KKP-133819, 147145, 147544, and 150151, which has been implemented with the support provided by the Ministry of Culture and Innovation of Hungary from the National Research, Development and Innovation Fund, financed under the ADVANCED-24 funding scheme. This research was funded by the grant 2024-1.2.8-TÉT-IPARI-CN-2025-00011,
with the support provided by the National Research,
Development and Innovation Office from the National Research,
Development and Innovation Fund, and financed under the
2024-1.2.8-TÉT-IPARI-CN funding scheme. 
}

\author{Florian Grundbacher}
\date{\today}

\begin{document}

\begin{abstract}
As a counterpart to classical vector balancing, we study large signed sums of unit vectors in finite-dimensional Minkowski spaces and develop their connection with polarization problems for convex bodies. 
For every convex body $K\subset\R^d$ containing the origin in its interior and every number $n\geq1$ of vectors, we show that the corresponding large signed sum and polarization constants coincide in a common quantity $\mu(K,n)$, which depends only on the symmetric core $K\cap(-K)$.

Our main result determines the sharp universal lower bound
\[
    \mu(K,n)
    \geq \mu(B_\infty^d,n)
    = \frac{1}{n}\left\lceil\frac{n}{d}\right\rceil,
\]
identifying parallelotopes as global minimizers for every $d$ and $n$. 
In the opposite direction, we prove general upper bounds that show, in particular, that the Euclidean ball is a maximizer up to an absolute constant factor. We also obtain sharp and nearly sharp results in several special cases, including the planar setting.

As $n\to\infty$, we prove that $\mu(K,n)$ converges to the Macphail constant of $K$, equivalently to its $1$-absolutely summing constant. Thus, the finite signed-sum problem provides a discrete counterpart of classical Banach space invariants. Combining this connection with sharp results on projection constants, we characterize equality in the corresponding upper bound for the Macphail constant in terms of maximal real equiangular tight frames.

Finally, our methods extend to arbitrary vector families, support functions of compact convex sets, and circumradii and diameters of Minkowski sums.
\end{abstract}

\maketitle

\section{Introduction}\label{sec:intro}

Vector balancing problems have been studied extensively for at least half a century (see, e.g., \cite{AB2024,barany2008Steinitz,dadush2018balancing}). In their classical form, one assigns signs to a given family of vectors so that the resulting signed sum is small in a prescribed norm. The complementary problem of forcing a signed sum to be large has received considerably less attention, except for the Euclidean setting \cite{AGM2021,AN2020,FuWangYan2023,Zhang2026Hidden}. In this paper, we study this large signed sum problem in finite-dimensional Minkowski spaces, i.e., $\R^d$ equipped with a Minkowski norm, and relate it to polarization problems for convex bodies, a subject with a substantial literature; see, e.g., \cite[Chapter~14]{borodachov2019discrete}. The corresponding connection for the Euclidean unit sphere has been established in \cite{AN2020,FuWangYan2023}.

We demonstrate that the resulting geometric quantities fit naturally into the framework of finite-dimensional Banach space theory.
For a convex body $K\subset\R^d$ and a fixed number $n$ of vectors, the large signed sum and polarization problems define a common constant $\mu(K,n)$.
We show that, as $n\to\infty$, these constants converge to the Macphail constant of the associated Minkowski space, equivalently to its $1$-absolutely summing constant.
Through the classical relation with projection constants, this places these finite vector problems within a well-developed Banach space framework.

Accordingly, several of our results may be viewed as precise, finite analogues of classical Banach space estimates in a more rigid setting.
We obtain sharp estimates that retain the dependence on both the dimension and the number of vectors, identify extremal configurations in several cases, and derive counterparts of classical inequalities.
In particular, we establish the precise, sharp universal lower bound for the polarization constants, and we prove general upper bounds, stability estimates, and extremal results connected with projection constants and equiangular tight frames.

\section{Terminology}\label{sec:terminology}

We start by introducing some basic geometric notions. We work with the class $\K_o^d$ of convex bodies (compact, convex subsets of $\R ^d$ with non-empty interior) that contain the origin $o$ in their interior.
For $X,Y \subset \R^d$, $t \in \R^d$, and $\nu \in \R$, we write $X+Y := \{ x + y : x \in X, y \in Y \}$ for the Minkowski sum, $t + X = \{t\} + X$ for the $t$-translation, and $\nu X := \{ \nu x : x \in X \}$ for the $\nu$-dilation. We abbreviate $-X := (-1) X$ and $X-Y := X + (-Y)$, and say that $X$ is {\em symmetric} if $-X = X$, i.e., $X$ is symmetric with respect to the origin.
The boundary and convex hull of $X$ are $\bd X$ and $\conv X$, respectively.
The closed segment connecting $x,y \in \R^d$ is denoted by $[x,y]$.
We write $[m] := \{1, \ldots, m\}$ for an integer $m \geq 1$.

For $p \geq 1$, the $\ell_p$-norm on $\R^d$ is denoted by $\|.\|_p$, while its unit ball is $B_p^d$. In the Euclidean ($p=2$) case, we simply write $|.| := \|.\|_2$ with associated inner product $\la ., . \ra$, and let $S^{d-1}:= \bd B_2^d$ denote the Euclidean unit sphere. For $K \in \K_o^d$, the Minkowski norm generated by $K$ (also called the $K$-norm) is
defined as 
\[
    \| x \|_K
    = \min \{\nu \geq 0: x \in \nu K \}
\]
for $x \in \R^d$.
Furthermore, the support function $h_K$ of $K$ is defined by
\[
h_K(x) = \max_{v \in K} \la x,v \ra
\]
for $x \in \R^d$. 
The polar $K^\circ$ of $K$ is, as usual, 
\begin{equation*}
    K^\circ
    = \{ x \in \R^d: \la x, v \ra \leq 1 \text{ for every } v\in K \}.
\end{equation*}
If $K \in \K_o^d$, then $K^\circ \in \K_o^d$ and $(K^\circ)^\circ = K$. Furthermore, it is well-known (cf.~\cite[Section~0.8]{gardner1995geometric}) that for any $K \in \K_o^d$ and $x \in \R^d$,
\begin{equation}\label{eq:polarsupport}
    h_{K^\circ}(x)
    = \| x\|_K.
\end{equation}
For more details, see \cite{gardner1995geometric} and \cite{schneider2013convex}.

Let us introduce the protagonists of the article.
The first central topic is the large signed sum problem, that is, maximizing the $K$-norm of signed sums of vectors on $\partial K$ (i.e.,~unit vectors in the $K$-norm). The Euclidean version
is discussed in \cite{AGM2021, AN2020, FuWangYan2023}.
This problem is dual to the classical notion of the signed sum constant (cf.~\cite[Section~4]{barany2008Steinitz}), where the goal is to find a signed sum of the given vectors whose norm is small. 

\begin{definition}\label{def:LSS}
For a convex body $K \in \K_o^d$ and an integer 
$n \geq 1$, the {\em $n$-th large signed sum constant of $K$} is given by
\begin{equation*}
    \mu(K, n)
    := \min_{u_1, \ldots, u_n \in \bd K} \ \max_{\eps \in \{ \pm 1 \}^n}  \frac{ \Big\| \sum_{i = 1}^n \eps_i u_i \Big\|_K}{n} ,
\end{equation*}
where $\eps=(\eps_1, \ldots, \eps_n) \in \{ \pm 1 \}^n$.
\end{definition}

\noindent
Large signed sum constants are also called {\em minimal energy} constants, see, e.g., \cite{FuWangYan2023, Zhang2026Hidden}.

We note that in the definition above, the vectors $u_1, \ldots, u_n$ need not be distinct. Accordingly, $\{ u_1, \ldots, u_n \}$ will be referred to as a {\em vector family}.

\medskip

The other central topic of the paper is the polarization problem for convex bodies.

\begin{definition}\label{def:pol}
For a convex body $K \in \K_o^d$ and an integer $n \geq 1$, the {\em $n$-th polarization constant of $K$} is given by 
\begin{equation}\label{eq:poldef}
    M(K,n)
    :=\min_{u_1, \ldots, u_n \in \bd K} \ \max_{v \in K^\circ } \frac {\sum_{i=1}^n |\la u_i, v \ra |}{n}.
\end{equation}
\end{definition}

Although it is not needed for the subsequent arguments, one can further generalize the above concept by introducing, for a convex body $K \in \K_o^d$, an integer $n \geq 1$, and for a scalar $p>0$, 
the {\em  $n$-th $\ell_p$-polarization constant of $K$} as 
\begin{equation*}
M_n^p(K) := \min_{u_1, \ldots, u_n \in \bd K} \ \max_{v \in K^\circ } \frac{\sum_{i=1}^n |\la u_i, v \ra | ^p}{n}.
\end{equation*}
Apart from the scaling factor, this is a special case of the two-plate Riesz polarization problem with exponent $(-p)$; see \cite[Section~14.2]{borodachov2019discrete}. 
Various notations have been used for analogous quantities -- in particular, $\mu(B_2^d,n)$ is denoted as $M_n^1(S^{d-1})$ in \cite{AN2020}.

As the compactness of $K$ and $K^\circ$ imply, 
$\mu(K,n)$ and $M(K,n)$ are well-defined, and they are bounded above by $1$ as a consequence of Lemma~\ref{lem:symm} and the triangle inequality. 
Note that both $n \cdot \mu(K,n)$ and $n \cdot M(K,n)$ are monotone increasing in $n$ and are invariant under invertible linear transformations applied to $K$. This also implies that $\mu(K,n) \geq \frac 1 n$ and $M(K,n) \geq \frac 1 n$ for every $n \geq 1$.

\section{Estimates for a Finite Number of Vectors}\label{sec:finite}

Our main objective is to give lower and upper bounds for the quantities $\mu(K,n)$ and $M(K,n)$.
The first result, which extends \cite[Proposition~3]{AN2020} to general Minkowski norms, states that these two constants are indeed equal.
While the observation is standard (see, e.g., \cite{AN2020,FuWangYan2023} for the Euclidean case and \cite{Garling_Gordon} for the general symmetric convex bodies), it is central to our approaches.
We therefore provide a short proof in Section~\ref{sec:reductions}.

\begin{theorem}\label{thm:muM}
For any 
$K \in \K_o^d$ and any $n \geq 1$, 
\[
    \mu(K,n)
    = M(K,n).
\]
\end{theorem}

In light of the statement above, from now on, we primarily use the notation $\mu(K,n)$, while referring to the quantity as the ($n$-th) polarization constant of $K$.

Next, we present two important reductions. The first implies
that in \eqref{eq:poldef}, the maximum  may be taken over the set of extremal points of $K^\circ$.

\begin{lemma}\label{lem:convhull}
For any $u_1, \ldots, u_n \in \R^d$ and any non-empty compact $V \subset \R^d$,
\begin{equation}\label{eq:polconvhull}
    \max_{v \in V} \sum_{i=1}^n |\la u_i, v \ra|
    = \max_{v \in \conv V} \sum_{i=1}^n |\la u_i, v \ra|.
\end{equation}
\end{lemma}
\noindent
This immediately simplifies the calculation of the polarization constant of polytopes, as one has to maximize only over the set of normal directions of facets.

We then show that the calculation of polarization constants may be reduced to symmetric convex bodies.
Therefore, we may assume in the forthcoming discussions that $K$ is symmetric.

\begin{lemma}\label{lem:symm}
For any  $K \in \K_o^d$ and any integer $n \geq 1$,
\[
    \mu(K,n)
    = \mu(K \cap (-K), n).
\]
\end{lemma}

After these reductions are shown in Section~\ref{sec:reductions}, we turn to stability estimates. 

\begin{lemma}\label{lemma:BM}
Suppose that $K,L \in \K_o^d$ with $L \subset TK \subset \nu L$ for some $\nu \geq 1$ and an invertible linear transformation $T$.
Then 
\[
    \max \Big\{ \frac{\mu(K,n)}{\mu(L,n)}, \frac{\mu(L,n)}{\mu(K,n)} \Big\}
    \leq \nu.
\]
\end{lemma}
In particular, if $K$ and $L$ are symmetric convex bodies in $\R^d$, then 
\begin{equation}\label{eq:BM}
    \max \Big\{ \frac{\mu(K,n)}{\mu(L,n)}, \frac{\mu(L,n)}{\mu(K,n)} \Big\}
    \leq d_{BM}(K,L),
\end{equation}
where $d_\mathrm{{BM}}(K,L)$ denotes the Banach--Mazur distance between $K$ and $L$ (cf.~\cite{artstein2015asymptotic}).
Therefore, the $n$-th polarization constant can be understood as a continuous mapping on the compact metric space of (classes of linear transformations of) symmetric convex bodies equipped with the (logarithm of the) Banach--Mazur distance. 
Consequently, it attains its minimum and maximum on this space, which are also its extrema on all of $\K^d_o$ by Lemma~\ref{lem:symm}.
We may therefore ask to determine the extremizers and the corresponding sharp bounds on the $n$-th polarization constant over $\K_o^d$.

Inequality \eqref{eq:BM} has two-fold applications in this context, as it leads to both lower and upper bounds. First, John's theorem on the approximation of convex bodies by ellipsoids \cite{john1948extremum} implies that any symmetric convex body in $\R^d$ has Banach--Mazur distance at most $\sqrt{d}$ from $B_2^d$ (see \cite{ball1997elementary}). Thus, by Lemma~\ref{lem:symm} and \eqref{eq:BM}, we derive that
for any $K \in \K_o^d$,
\begin{equation}\label{eq:pol_lower_ball}
    \mu(K, n)
    \geq \frac 1 {\sqrt{d}} \cdot \mu(B_2^d,n).
\end{equation}
Therefore, the already established estimates for the polarization constants of $B_2^d$ lead to general lower bounds. These constants are known precisely in special cases: for small values of $n$,
and in the planar case $d=2$. First, for $n \leq d$,
it has been shown in \cite[Section~4]{AN2020} (see also \cite[Proposition~2.2]{FuWangYan2023}) that
\begin{equation}\label{eq:euclidean_small_n}
    \mu(B_2^d,n)
    = \frac 1 {\sqrt{n}}.
\end{equation}
Second, it has been proven in  {\cite[Proposition~5]{AN2020}}, and independently in \cite[Proposition~2.4]{FuWangYan2023}, that the minimum in the definition of $\mu(B_2^2, n)$ corresponds to sets of $n$ vectors that are equally distributed on a half-circle. We provide a new, direct proof for this statement in Section~\ref{sec:Euclidean}.

\begin{proposition}[\cite{AN2020,FuWangYan2023}]\label{prop:planar}
For any $n \geq 1$,
\[
    \mu(B_2^2, n)
    = \frac{1}{n \sin(\frac \pi {2 n})}.
\]
The minimum is achieved precisely by vectors $u_1, \ldots, u_n \in S^1$ for which $\{\pm u_1, \ldots, \pm u_n\}$ forms the vertex set of a regular $2n$-gon.
\end{proposition}

In the general case, we know the following estimates for $\mu(B_2^d,n)$ from \cite{AN2020}. We use the standard notation $\kappa_d$ for the volume of $B_2^d$:
\[
    \kappa_d
    = \frac{\pi^{d/2}}{\Gamma(\frac d 2 + 1)}.
\]

\begin{proposition}[\cite{AN2020}]\label{prop:polball}
For any $d \geq 1$ and any $n \geq d$,   
\begin{equation}\label{eq:ballpol}
    \frac{2 \, \kappa_{d-1}}{\sqrt{d} \, \kappa_d} \cdot \frac 1 {\sqrt{d}}
    \leq \mu(B_2^d,n)
    \leq \frac{1}{\sqrt{d}}.
\end{equation}
\end{proposition}
\noindent
A short proof for the proposition is presented in Section~\ref{sec:Euclidean}. The asymptotic behavior of the lower bound is explained by
\begin{equation}\label{eq:ballconstant}
    \frac{2 \, \kappa_{d-1}}{\sqrt{d} \, \kappa_d}
    = \sqrt{ \frac 2 \pi } + O \bigg( \frac{1}{d} \bigg)
\end{equation}
(here and later on, asymptotic orders are understood in terms of $d \to \infty$).
Moreover, the above quantity is decreasing in $d$. We note that asymptotically matching estimates are also proved in \cite[Theorem~3.3~\&~Theorem~3.4]{FuWangYan2023}.

We also remark that both bounds in \eqref{eq:ballpol} are sharp:
while the lower bound is attained only if $d = 1$,
as $n \to \infty$ for a fixed $d \geq 2$,
sets of nearly evenly distributed points on a half-sphere lead to values that behave asymptotically like 
the lower estimate in \eqref{eq:ballpol}.
For $n = d$, the upper bound is attained.

Now, the combination of \eqref{eq:pol_lower_ball} and Proposition ~\ref{prop:polball} leads to the asymptotic lower bound 
\begin{equation*}
    \mu(K,n)
    \geq \Omega \bigg( \frac{1}{d} \bigg)
\end{equation*}
that holds for every $K \in \K_o^d$.
We make this estimate precise in Section~\ref{sec:lower}
by showing that the cube is a global minimizer of the polarization constant for any dimension and
any fixed number of vectors.

\begin{theorem}\label{thm:sharplower}
For any $K \in \K_o^d$ and any $n \geq 1$, 
\[
    \mu(K,n)
    \geq \mu(B_\infty^d, n)
    = \frac 1 n \cdot \bigg\lceil \frac n d \bigg\rceil.
\]
\end{theorem}

As a natural counterpart, we turn to upper bounds. As we noted before,  Lemma~\ref{lemma:BM} can also be used to provide upper estimates. In particular, a result of Giannopoulos~\cite{giannopoulos1995note} states that any symmetric convex body in $\K^d_o$ has distance at most $O(d^{5/6})$ from the cube $B_\infty^d$.
This, via \eqref{eq:BM} together with Lemma~\ref{lem:symm} and Theorem~\ref{thm:sharplower}, 
implies that
for any
$K \in  \K_o^d$, 
\[
    \mu(K,n)
    \leq O \bigg( \frac{1}{d^{1/6}} \bigg).
\]
We improve on this estimate substantially in Section~\ref{sec:upper}:

\begin{theorem}\label{thm:upperbd}
For any $K \in \K_o^d$ and any $n \geq 1$, 
\begin{equation}\label{eq:upperbd}
    \mu(K,n)
    \leq  2 \max \Big\{ \frac 1 {\sqrt{n}}, \frac{1} {\sqrt{d}} \Big\}.
\end{equation}
Moreover,
\begin{equation}\label{eq:upperbd2}
    \mu(K,n)
    < \frac{1}{\sqrt{d}} + \frac{d (d+1)}{2n}.
\end{equation}
\end{theorem}

The first estimate is stronger when $n$ is small relative to $d$, whereas the second bound is smaller when $d$ is fixed and $n$ is sufficiently large.

Combining the first estimate \eqref{eq:upperbd} with \eqref{eq:ballpol} and \eqref{eq:ballconstant} leads to the bound
\[
    \mu(K, n)
    < \sqrt{2 \pi} \cdot \mu(B_2^d, n),
\]
which shows that the Euclidean ball is an asymptotic maximizer. It remains an open question to determine the precise maxima:

\begin{problem}\label{problem1}
For any given $d, n \geq 2$, determine the maximum of $\mu(K, n)$ over $K \in \K_o^d$.
\end{problem}

We solve special cases of Problem~\ref{problem1} in Section~\ref{sec:special}.
The first one is essentially covered by \cite[Remark]{dvoretzky_rogers},
for which we provide a new, slightly more direct proof.

\begin{proposition}[\cite{dvoretzky_rogers}]\label{prop:n=2}
For any $K \in \K_o^d$,
\[
    \mu(K,2)
    \leq \mu(B_2^d,2)
    = \begin{cases}
        1, & \text{if $d=1$,} \\
        \frac{1}{\sqrt{2}}, & \text{else.}
    \end{cases}
\]
\end{proposition}

\begin{proposition}\label{prop:planar_n=3}
For any $K \in \K_o^2$,
\[
    \mu(K,3)
    = \frac 2 3.
\]
\end{proposition}

Combining the above two propositions, we arrive at (almost) sharp upper bounds for general $n \geq 1$ in the planar case.

\begin{corollary}\label{cor:planar_upper}
For any $K \in \K_o^2$ and any $n \geq 1$,
\[
    \mu(K,n)
    \leq \begin{cases}
        1, & \text{if } n=1,
        \\
        \frac{2}{3}, & \text{if $n \geq 2$ and $n \equiv 0$ (mod $3$),}
        \\
        \frac{2}{3} + \frac{2}{n} \big( \sqrt{2} - \frac{4}{3} \big), & \text{if $n \geq 2$ and $n \equiv 1$ (mod $3$),}
        \\
        \frac{2}{3} + \frac{1}{n} \big( \sqrt{2} - \frac{4}{3} \big), & \text{if $n \geq 2$ and $n \equiv 2$ (mod $3$).}
    \end{cases}
\]
The inequality is sharp when $n \leq 3$ or $n \equiv 0$ (mod $3$).
\end{corollary}

\section{Connection with Banach Space Estimates}\label{sec:Banach}

It is natural to consider the limits of the quantities $\mu(K,n)$ and $M(K,n) $ as the number of vectors $n$ tends to infinity. As we shall see, one may also relax the normalization condition $u_i \in \partial K$ for every $i$ at the same time. Starting from the large signed sum problem (cf.~Definition~\ref{def:LSS}), this leads to the notion of {\em Macphail constants}, which we introduce via the Kadets--Snobar variant~\cite{Kadets_Snobar}.

\begin{definition}\label{def:macphail}
For a convex body $K \in \K_o^d$, the {\em Macphail constant} is given by
\begin{equation*}
    \mu(K)
    := \inf_{ X \subset \R^d \setminus \{ o \} \text{ finite}} \ \max_{ \eps \in \{ \pm 1 \}^X} \frac{\| \sum_{x \in X} \eps(x) \cdot x \|_K}{ \sum{_{x \in X}} \| x\|_K }.
\end{equation*}
\end{definition}
\noindent
Above, the infimum could also be taken over all finite vector {\em families} (i.e., with possible repetitions),
as replacing $k$ copies of a vector $x$ by a single vector $k \cdot x$ does not change the value on the right-hand side.

We note that for symmetric $K$, Macphail constants are frequently defined alternatively by 
\[
    \inf_{ X \subset \R^d \setminus \{ o \} \text{ finite}} \ \max_{Y \subset X} \frac{\| \sum_{y \in Y}  y \|_K}{ \sum{_{x \in X}} \| x\|_K },
\]
see, e.g., \cite{Gordon68}.
The above quantity equals $\frac 1 2 \cdot \mu(K)$ in our terminology
(see \cite[Proposition]{Garling_Gordon}).

On the other hand, the continuous analogue of the polarization problem (cf.~Definition~\ref{def:pol}) leads to the concept of {\em 1-absolutely summing constants:}
\begin{definition}
For a convex body $K \in \K_o^d$, its {\em $1$-absolutely summing constant} is given by
\[
    \mu_1(K)
    := \inf_{ X \subset \R^d \setminus \{ o \} \text{ finite}} \ \max_{v \in K^\circ} \frac{\sum_{x \in X} |\la x, v \ra|}{ \sum{_{x \in X}} \| x\|_K }.
\]
\end{definition}

We prove in Section~\ref{sec:Limit} that, as one would naturally expect, the large signed sum and polarization constants converge to the Macphail and $1$-absolutely summing constants, respectively. Accordingly, the former may be viewed as finite, more rigid counterparts of these limiting quantities.

\begin{theorem}\label{thm:pollim}
For any $K \in \K_o^d$, 
\[
    \lim_{n \to \infty} \mu(K, n)
    = \inf_{n \geq 1} \mu(K, n)
    = \mu (K)
\]
and 
\[
    \lim_{n \to \infty} M(K, n)
    = \inf_{n \geq 1} M(K, n)
    = \mu_1 (K).
\]
\end{theorem}
\noindent Note that even though the limit and the infimum of $\mu(K, n)$ as $n \to \infty$ coincide, the sequence is not always decreasing in $n$, as Theorem~\ref{thm:sharplower} shows. Also,  Theorem~\ref{thm:muM} immediately implies that $\mu_1(K) = \mu(K)$ for every $K \in \K_o^d$. This identity was established previously for symmetric convex bodies; see~\cite{Garling_Gordon, Gordon69}.

Theorem~\ref{thm:pollim} establishes a connection between $\mu(K,n)$ and $\mu(K)$ that yields estimates in both directions. On the one hand, knowledge of the precise value of $\mu(K)$ provides a lower bound for  $\mu(K,n)$ for every $n \geq 1$, along with an asymptotic upper bound as $n \to \infty$. Several such examples are discussed in~\cite{Gordon69}, where the $1$-absolutely summing constants are determined for a number of classical examples. 
Furthermore, since Definition~\ref{def:macphail} imposes no normalization condition on the vectors, this approach also extends the preceding results to families of vectors that are not necessarily of unit norm. For example, combining Proposition~\ref{prop:polball} with Definition~\ref{def:macphail} immediately leads to the following estimate. 

\begin{corollary}\label{cor:eucl_general_vectors}
For any $u_1, \ldots, u_n \in \R^d$,
\[
    \max_{\varepsilon \in \{\pm1\}^n} \Big| \sum_{i=1}^n \varepsilon_i u_i \Big|
    \geq \frac{2 \, \kappa_{d-1}}{d \, \kappa_d} \sum_{i=1}^n |u_i|,
\]
with equality if and only if $d=1$ or all $u_i$ are zero.
Moreover, the inequality is sharp for any $d \geq 1$ as $n \to \infty$.
\end{corollary}

Conversely, letting $n \to \infty$ in our preceding results, we obtain information about $\mu(K)$. 
As the first example, we note that the Macphail and $1$-absolutely summing constants were originally defined in the setting of normed spaces and, hence, for symmetric convex bodies. However, Lemma~\ref{lem:symm} shows that this symmetry assumption entails no genuine loss of generality, since
\begin{equation}\label{eq:mu_symm}
    \mu(K)
    = \mu(K \cap -K)
\end{equation}
for every $K \in \K_o^d$.

Second, Theorem~\ref{thm:sharplower} implies the sharp lower bound $ \mu(K)\geq \frac 1 d$ for the Macphail and $1$-absolutely summing constants of symmetric convex bodies, which was proved by Kadets and Snobar~\cite{Kadets_Snobar}; see also Gordon~\cite{Gordon69}. Among symmetric convex bodies, the equality case was characterized by Deschaseaux~\cite{Deschaseaux1973}; see also Garling~\cite{Garling1974Trace}. The extension to non-symmetric convex bodies follows from \eqref{eq:mu_symm}.

\begin{proposition}
For any convex body $K \in \K_o^d$,
\begin{equation}\label{eq:muest_lower}
    \mu(K)
    \geq \frac{1}{d},
\end{equation}
with equality
if and only if $K \cap (-K)$ is a parallelotope, i.e., a linear image of $B_\infty^d$.
\end{proposition}

We note that, to the best of our knowledge, no elementary proof of the characterization of the equality cases is available in the literature.

Turning to upper estimates, letting $n\to\infty$ in \eqref{eq:upperbd2} of Theorem~\ref{thm:upperbd} yields
\[
    \mu(K)
    \leq \frac 1 {\sqrt{d}}.
\]
This agrees asymptotically with the following classical estimate of König and Tomczak-Jaegermann~\cite[Proposition~3]{KonigTomczak1990}. Below, and later on, 
\begin{equation}\label{eq:gammad}
    \gamma_d
    := \frac{2 + (d-1)\sqrt{d+2}}{d(d+1)}.
\end{equation}

\begin{proposition}[\cite{KonigTomczak1990}]\label{prop:uppermu1}
For any convex body $K \in \K_o^d$,
\[
    \mu(K)
    \leq \gamma_d .
\]
\end{proposition}

Similarly to the lower bound, our next goal is to characterize the cases of equality above.
To that end, we introduce the notion of  {\em projection constants}, another well-studied concept describing Banach spaces~\cite{grunbaum1960projection}:
\begin{definition}
Let $K\in\K_o^d$ be a symmetric convex body, and let
$X=(\R^d,\|\cdot\|_K)$ be the real $d$-dimensional Banach space
whose unit ball is $K$. The \emph{projection constant} of $K$ is defined by
\begin{align*}
    \lambda(K)
    := \inf\{ \lambda > 0 :
        & \text{ from every Banach space $Y \supset X$, there exists} \\ & \text{ a surjective linear projection onto $X$ with norm at most $\lambda$}
    \}.
\end{align*}
\end{definition}
The following connection of the projection constant with the 1-summing and Macphail constants was established by Gordon.
\begin{proposition}[{\cite{Gordon68}}]\label{prop:gordon}
For any symmetric $K \in \K_o^d$,
\[
    \mu(K)
    \leq \frac{\lambda(K)}{d}.
\]
\end{proposition}

Thus, upper bounds on $\lambda(K)$ yield corresponding estimates for $\mu(K)$. The relationship between the five main quantities studied in the paper is summarized in the diagram below. 
\[
\begin{array}{c@{\,\,}c@{\,\,}c@{\;\;}c}
    \mu(K,n) & \xrightarrow{\ n\to\infty\ } & \mu(K) &
    \\
    [-0.1em] \rotatebox{90}{$=$} & & \rotatebox{90}{$=$} & \displaystyle \leq \ \frac{\lambda(K)}{d}
    \\
    [-0.1em] M(K,n) & \xrightarrow{\ n\to\infty\ } & \mu_1(K) &
\end{array}
\]

In the planar case, the inequality $\lambda(K) \leq \frac{4}{3}$, which is a special case of Proposition~\ref{prop:proj_strong}, was conjectured by Grünbaum~\cite{grunbaum1960projection}. It was proved only fifty years later by Chalmers and Lewicki~\cite{chalmers2010}; see also the simplified proofs in~\cite{basso2019computation,DL23}. The corresponding higher-dimensional estimate, which matches Proposition~\ref{prop:uppermu1}, was first stated in~\cite{KonigTomczak1994} and verified recently by Der\c{e}gowska and Lewandowska~\cite{DL23}, while the equality cases were characterized by Kobos~\cite[Theorem~1.2]{Kobos25}.

\begin{proposition}[\cite{DL23, Kobos25}]\label{prop:proj_strong}
For any symmetric $K \in \K_o^d$,
\[
    \lambda(K)
    \leq d \gamma_d,
\]
with equality if and only if there exist a maximal equiangular tight frame $w_1, \ldots, w_N \in \R^d$ with $N = \frac{d(d+1)}{2}$, and an invertible linear transformation $T$ such that 
\begin{equation*}
    \conv \{\pm \etfvar_1,\ldots, \pm \etfvar_\etfind\} \subset T( K^\circ )
    \subset \frac 1 {N \gamma_{d}} \sum_{i=1}^{\etfind} [-\etfvar_i,\etfvar_i].
\end{equation*}
\end{proposition}

\noindent
Vectors $\etfvar_1, \ldots, \etfvar_\etfind \in \R^d$ form a {\em tight frame} if there exists $\alpha > 0$ with
\begin{equation}\label{eq:TF1}
    \sum_{i=1}^\etfind \la x, \etfvar_i \ra \etfvar_i
    = \alpha x
\end{equation}
for all $x \in \R^d$. Equivalently, 
\begin{equation}\label{eq:TF}
    \sum_{i=1}^\etfind \etfvar_i \otimes \etfvar_i
    = \alpha \, I_d
\end{equation}
where $ \etfvar \otimes  \etfvar =  \etfvar  \etfvar^\top$ is the rank $1$ operator on $\R^d$ associated with $ \etfvar \in \R^d$ and $I_d$ is the identity operator on $\R^d$. In case $| \etfvar|=1$, $ \etfvar \otimes  \etfvar$ is the orthogonal projection onto the linear hull of $w$. 
The tight frame is {\em equiangular} (ETF)
if all $\etfvar_i$ have the same length and there exists some $\beta \geq 0$ with
\begin{equation*}
    |\la \etfvar_i, \etfvar_j \ra|
    = \beta
\end{equation*}
for all $i,j \in [\etfind]$, $i \neq j$. The cardinality of an ETF is guarded from below by comparing the ranks in \eqref{eq:TF}, and from above by Gerzon's bound (cf.~\cite{LemmensSeidel1973}):
\begin{equation}\label{eq:Nmax}
    d \leq \etfind
    \leq \frac{d (d+1)}{2}.
\end{equation}
The ETF is said to be {\em maximal} if the upper bound is reached.
The only dimensions in which maximal real ETFs are currently known to exist
are $d=1,2,3,7,23$, and evidence suggests that this list may be complete
\cite{BannaiMunemasaVenkov2005,Kobos25}.

We are ready to characterize the cases of equality  of the upper bound on Macphail constants stated in Proposition~\ref{prop:uppermu1}. The extremal examples are given by  polar bodies associated with ETFs that also yield the sharpness of the estimate of Proposition~\ref{prop:proj_strong}, see~\cite[Theorem~1.2]{Kobos25}. Generalizing \eqref{eq:gammad}, we introduce the quantity
\begin{equation*}
    \gamma_{d,\etfind}
    := 1+\sqrt{\frac{(\etfind-d)(\etfind-1)}{d}} . 
\end{equation*}
Comparing this with \eqref{eq:gammad} shows that for $N= \frac{d(d+1)}{2}$, we have $\gamma_{d,\etfind} = N \gamma_d$.

\begin{proposition}\label{prop:etf}
Suppose that the unit vectors $\etfvar_1, \ldots, \etfvar_\etfind \in S^{d-1}$ form an equiangular tight frame in $\R^d$ and that $K \in \K_o^d$ satisfies
\begin{equation}\label{eq:Kpolcond}
    \conv \{\pm \etfvar_1,\ldots, \pm \etfvar_\etfind\}
    \subset T( (K \cap -K)^\circ )
    \subset \frac 1 {\gamma_{d,\etfind}} \sum_{i=1}^{\etfind} [-\etfvar_i,\etfvar_i]
\end{equation} 
for some invertible linear transformation $T$.
Then, for any $n \geq 1$,
\begin{equation}\label{eq:muKn_lower}
    \mu(K,n)
    \geq \frac{\gamma_{d,N}}{\etfind},
\end{equation}
with equality when equality holds on the left-hand side in \eqref{eq:Kpolcond} and $n$ is a multiple of $N$.
\end{proposition}

\noindent
By taking polars, the condition \eqref{eq:Kpolcond} is equivalent to 
\[
    \Big \{x \in \R^d: \sum_{i=1}^N |\la x, w_i \ra| \leq  \gamma_{d,N}   \Big \}  \subset \widetilde{T} (K \cap (-K))
    \subset \Big \{ x \in \R^d: |\la x, w_i \ra | \leq 1 \text{ for all } i \in [N] \Big \}.
\]
for some invertible linear transformation $\widetilde{T}$.
Let us also point out that both assumptions for the equality case are in general needed:
When $N=d$, the ETF is simply an orthonormal basis, so that \eqref{eq:Kpolcond} is always true by the Auerbach lemma. Moreover, the polar of the left-hand set in \eqref{eq:Kpolcond} is then just a cube, for which we know from Theorem~\ref{thm:sharplower} that its $n$-th polarization constant is $\frac{1}{d}$ only if $n$ is a multiple of $d = N$.

If the maximum cardinality \eqref{eq:Nmax} of an ETF is reached, the lower bound in \eqref{eq:muKn_lower} agrees with the upper bound in Proposition~\ref{prop:uppermu1}. Combined with Proposition~\ref{prop:gordon} and Proposition~\ref{prop:proj_strong}, this yields the following sharp upper bound on the Macphail constant of a convex body for this specific set of dimensions, along with the characterization of equality cases.

\begin{theorem}\label{thm:mu_upper}
For any $K \in \K_o^d$, 
\begin{equation}\label{eq:mu_uppergamma}
    \mu (K)
    \leq \gamma_d,
\end{equation}
with equality if and only if there exist a maximal equiangular tight frame $w_1, \ldots, w_N \in \R^d$ with $N = \frac{d(d+1)}{2}$, and an invertible linear transformation $T$ such that 
\begin{equation}\label{eq:mumax_characterization}
    \conv \{\pm \etfvar_1,\ldots, \pm \etfvar_\etfind\} \subset T( (K \cap -K)^\circ )
    \subset \frac 1 {N \gamma_{d}} \sum_{i=1}^{\etfind} [-\etfvar_i,\etfvar_i].
\end{equation}
\end{theorem}

\noindent In particular, for $d=2$, the bounding convex bodies in \eqref{eq:mumax_characterization} are equal: both are regular hexagons with edge length  $2$, whose polar is again a regular hexagon. Consequently, equality holds if and only if $K \cap (-K)$ is a symmetric, affine-regular hexagon. However, for larger dimensions, the convex hull on the left side of \eqref{eq:mumax_characterization} differs from the zonotope on its right side. Therefore, equality is attained for a wider class of convex bodies. For instance, in the $3$-dimensional case, a symmetric $K$ maximizes the  Macphail constant if and only if a linear transformation of $K$ fits between a regular icosidodecahedron and the regular dodecahedron whose facets are supersets of the pentagonal facets of the icosidodecahedron.

For dimensions $d =  7, 23$, the symmetric maximizers of the Macphail constant include the  polar of the Gosset polytope $3_{21}$, and the polar of the ``$Co_3$-polytope'', respectively.
The latter is the polytope in $\R^{23}$ whose vertices are the Euclidean unit vectors lying on the $276$ equiangular lines in $\R^{23}$ derived from the Leech lattice in $\R^{24}$ and whose symmetry group is given by the Conway group $Co_3$ (see, e.g., \cite{GoethalsSeidel1975}).

The continuous analogue of Problem~\ref{problem1}, that is, the question of exact maximizers of the Macphail constant, remains open in dimensions $d \neq 1,2,3,7,23$:
\begin{problem}\label{problem2}
For any given $d \geq 2$, determine the maximum of $\mu(K)$ over $K \in \K^d_o$.
\end{problem}

We conclude the article in Section~\ref{sec:extensions} with several extensions of the above results. By slight modifications of the proofs, we establish sharp generalizations for general vector sets, compact convex sets not containing the origin, and circumradii of Minkowski sums.

\section{Reductions and Relations}
\label{sec:reductions}

First, we establish the connection between the polarization and large signed sum constants.
It is largely a consequence of \eqref{eq:polarsupport}.

\begin{proof}[Proof of Theorem~\ref{thm:muM}]
To begin with, note that for any $w \in \R^d$, \eqref{eq:polarsupport} states that
\begin{equation}\label{eq:wv}
    \max_{v \in K^\circ} \la w, v \ra
    = \|w\|_K.
\end{equation}
Let now $ u_1, \ldots, u_n \in \bd K$ be arbitrary. For any $\eps \in \{ \pm 1 \}^n$, \eqref{eq:wv} leads to
\begin{equation}\label{eq:skmk}
    \Big\| \sum_{i = 1}^n \eps_i u_i \Big\|_K
    = \max_{v \in K^\circ} \Big \la \sum_{i = 1}^n \eps_i u_i , v \Big \ra
    = \max_{v \in K^\circ} \sum_{i = 1}^n \la \eps_i u_i, v \ra  \leq \max_{v \in K^\circ} \sum_{i = 1}^n | \la u_i, v \ra |.
\end{equation}

On the other hand, let $v \in K^\circ$ be a maximizer of $ \sum_{i = 1}^n | \la u_i, v \ra |$, and for each $i \in [n]$ 
set $\eps_i \in \{\pm 1\}$ such that $\la \eps_i u_i, v \ra \geq 0$. Then equality holds in \eqref{eq:skmk}. Combining the above facts, we have that
\[
\max_{\eps \in \{ \pm 1 \}^n} \Big\| \sum_{i = 1}^n \eps_i u_i \Big\|_K = \max_{v \in K^\circ} \sum_{i = 1}^n | \la u_i, v \ra |.
\]
Since the choice of $u_1, \ldots, u_n \in \bd K$ was arbitrary, this shows that 
\[
\mu(K,n) = M(K,n). 
\qedhere
\]
\end{proof}

We continue with the reduction steps: Lemmas~\ref{lem:convhull}~and~\ref{lem:symm}.
For the latter, we use the more detailed insight into the contributions of the individual vectors provided by the definition of the polarization constants.

\begin{proof}[Proof of Lemma~\ref{lem:convhull}]
By containment, it is straightforward that the left-hand side of \eqref{eq:polconvhull} is not larger than its right-hand side. Thus, it suffices to show that
\begin{equation}\label{eq:convhull}
    \max_{v \in V} \sum_{i=1}^n |\la u_i, v \ra|
    \geq \max_{w \in \conv V} \sum_{i=1}^n |\la u_i, w \ra|.
\end{equation}
To do so, let $w$ be a maximizer of the right-hand side,
which exists by compactness of $\conv V$ from compactness of $V$.
Since $w \in \conv V$, there exist some $v_1, \ldots, v_m \in V$ and convex coefficients $\lambda_1, \ldots, \lambda_m \geq 0$, $\sum_{j=1}^m\lambda_j =1$, for which $w = \sum_{j=1}^m \lambda_j v_j$. Then 
\[
    \sum_{i=1}^n |\la u_i, w \ra|
    = \sum_{i=1}^n \Big| \la u_i, \sum_{j=1}^m \lambda_j v_j \ra \Big|
    = \sum_{i=1}^n \Big| \sum_{j=1}^m \lambda_j \la u_i, v_j \ra \Big|
    \leq \sum_{j=1}^m \lambda_j \sum_{i=1}^n |\la u_i, v_j \ra|.
\]
Thus, 
$\sum_{i=1}^n |\la u_i, v_j \ra| \geq \sum_{i=1}^n |\la u_i, w \ra|$ for some $j \in [m]$, which implies \eqref{eq:convhull}.
\end{proof}

\begin{proof}[Proof of Lemma~\ref{lem:symm}]
Note that $(K \cap (-K))^\circ = \conv \{ K^\circ \cup (- K^\circ) \}$. Furthermore, observe that for every point $w \in \bd (K \cap (-K))$, there exists $\varepsilon \in \{\pm1\}$ such that $\varepsilon w \in \bd K$.
Accordingly, Lemma~\ref{lem:convhull} implies that 
\begin{align*}
    \mu(K \cap (-K), n)
    &= \min_{w_1, \ldots, w_n \in \bd (K \cap (-K))} \ \max_{z \in (K \cap (-K))^\circ } \frac 1 n \sum_{i=1}^n |\la w_i, z \ra |\\
    &=  \min_{w_1, \ldots, w_n \in \bd (K \cap (-K))} \ \max_{v \in K^\circ } \frac 1 n \sum_{i=1}^n |\la w_i, v \ra | \\
    &=  \min_{u_1, \ldots, u_n \in ( \bd K \cap \bd (K \cap (-K)))} \ \max_{v \in K^\circ } \frac 1 n \sum_{i=1}^n |\la u_i, v \ra |.
\end{align*}
On the other hand, in \eqref{eq:poldef}, if for some $u_i \in \bd K$ we have that $- u_i \notin K$, then replacing $u_i$ by $-u_i / \| - u_i \|_K$ clearly does not increase the maximum. Hence, one may restrict $\bd K$ to $\bd K \cap \bd (K \cap (-K))$ when minimizing over $n$-tuples for $\mu(K,n)$. Therefore, the quantity above equals $\mu(K,n)$, and the proof is complete.
\end{proof}

Finally, we prove the stability estimate, where we again rely on the definition of polarization constants to have better control over the contributions of the individual vectors. 

\begin{proof}[Proof of Lemma~\ref{lemma:BM}] Since $\mu(TK, n) = \mu(K, n)$, we may assume that  $L \subset K \subset \nu L$. By polarity, this  implies that 
\[
\frac 1 \nu L^\circ \subset K^\circ \subset L^\circ.
\]
Accordingly, for any set of vectors $u_1, \ldots, u_n \in \bd K$, 
\begin{equation*}
    \max_{v \in K^\circ} \sum_{i=1}^n |\la u_i, v \ra|
    \geq \max_{v \in K^\circ} \sum_{i=1}^n \Big| \Big \la \frac{u_i}{\| u_i \|_L}, v \Big \ra \Big|
    \geq \frac 1 \nu \cdot \max_{w \in L^\circ} \sum_{i=1}^n \Big| \Big \la \frac{u_i}{\| u_i \|_L}, w \Big \ra \Big|
    \geq \frac n \nu \cdot \mu(L, n).
\end{equation*}
This shows $\nu \cdot \mu(K,n) \geq \mu(L,n)$. Since also $K \subset \nu L \subset \nu K$, we may reverse the roles of $K$ and $L$ to obtain $\nu \cdot \mu(L,n) \geq \mu(K,n)$ and complete the proof.
\end{proof}

\section{Lower Bounds}
\label{sec:lower}

In this section, we prove that parallelotopes minimize the $n$-th polarization constant.
Our key insight is a structural observation about well-chosen generating systems that follows from an Auerbach-type argument.

\begin{lemma}\label{lemma:Auerbach}
For any set of vectors $u_1, \ldots, u_m \in \R^d$ with $m \geq d$, there exist some pairwise distinct indices $i_1, \ldots, i_d$ with
\begin{equation}\label{eq:good_basis}
    u_1, \ldots, u_m
    \in \sum_{j=1}^d [-u_{i_j},u_{i_j}].
\end{equation}
\end{lemma}

\begin{proof} Let $U$ be the linear subspace spanned by the $u_i$ and write $d' := \dim(U)$.
Choose $i_1, \ldots, i_{d'}$ from $[m]$ so that $T := \conv \{ o, u_{i_1}, \ldots, u_{i_{d'}} \}$ has maximum possible $d'$-dimensional volume among $d'$-dimensional simplices with one vertex at the origin $o$ and $d'$ vertices among $u_1, \ldots, u_m$. Observe that due to maximality, 
\[
    \conv \{\pm u_1, \ldots, \pm u_m\}
    \subset \sum_{j=1}^{d'} [-u_{i_j}, u_{i_j}]
\]
(see \cite{AHAK} for a detailed proof). In case $d' < d$, we add $\sum_{j=d'+1}^d [-u_{i_j},u_{i_j}] \ni o$ for any further $i_{d'+1}, \ldots, i_d$. This can only increase the right-hand set, which yields \eqref{eq:good_basis}.
\end{proof}

\begin{proof}[Proof of Theorem~\ref{thm:sharplower}]
We start with $\mu(B_\infty^d,n) \leq \frac 1 n \cdot \lceil \frac{n}{d} \rceil$.
Let $\{e_j\}_{j=1}^d$ be the standard orthonormal basis of $\R^d$. 
Then, for any $\varepsilon \in \{\pm1\}^d$,
\[
    \Big\| \sum_{j=1}^d \varepsilon_j e_j \Big\|_{\infty}
    = 1.
\]
Now, set $r := \lceil \frac{n}{d} \rceil$, $m := r \cdot d$, and let $u_1, \ldots, u_m$ be the vector system that consists of $r$ copies of $e_1, \ldots, e_d$.
Since $u_1, \ldots, u_m \in \bd B_\infty^d$, and $ n \cdot \mu(B_\infty^d,n) \leq m \cdot\mu(B_\infty^d,m)$ by monotonicity, we obtain
\[
    n \cdot  \mu(B_\infty^d,n)
    \leq m \cdot \mu(B_\infty^d,m)
    \leq  \max_{\varepsilon \in \{\pm1\}^m} \Big\| \sum_{i=1}^m \varepsilon_i u_i \Big\|_{\infty}
    \leq  r .
\]
It follows
$\mu(B_\infty^d,n) \leq \frac 1 n \cdot \lceil \frac{n}{d} \rceil$ as claimed.

Next, we show $n \cdot \mu(K,n) \geq r :=\lceil \frac{n}{d} \rceil$
for general $K \in \K^d_o$. Note that this holds trivially for $n \leq d$ 
by the monotonicity of $n \cdot \mu(K,n)$ in $n$. Hence, we assume $r \geq 2$ from now on.

Starting from any set of $u_1, \ldots, u_n \in \bd K$, 
we repeatedly select and drop $d$ vectors according to Lemma~\ref{lemma:Auerbach} so that in each step, the parallelotope generated by the selected $d$ vectors contains all the remaining vectors. By reindexing, we may thus assume that 
\[
    u_j \in \sum_{i=k d + 1}^{(k+1) d} [-u_i,u_i]
\]
for all $0 \leq k \leq r-2$ and $j > (k+1) d$. Note that in particular, containment holds in each step for $j = n$. Furthermore, 
\[
    u_n
    \in \sum_{i=(r-1) d + 1}^{n} [-u_i,u_i].
\]
Therefore,
\[
    r u_n
    \in \sum_{i=1}^n [-u_i,u_i]
    = \conv \Big\{ \sum_{i=1}^n \eps_i u_i : \varepsilon \in \{\pm1\}^n \Big\}.
\]
Since $K$ is convex, the polytope $\conv \big\{ \sum_{i=1}^n \eps_i u_i : \varepsilon \in \{\pm1\}^n \big\}$ must also have a vertex of large $K$-norm:
\begin{equation}\label{eq:convexity_ineq}
    \max_{\varepsilon \in \{\pm1\}^n} \Big\| \sum_{i=1}^n \varepsilon_i u_i \Big\|_K
    \geq \| r u_n \|_K
    = r.
\end{equation}
Accordingly,
\[
    \mu(K,n)
    \geq \frac r n
    = \frac 1 n \cdot \bigg\lceil \frac{n}{d} \bigg\rceil. \qedhere
\]
\end{proof}

\begin{remark}\label{rem:asymptotic_lower}
Once Theorem~\ref{thm:pollim} is established, the slightly weaker lower bound $\frac 1 d$ follows immediately from the Kadets--Snobar estimate~\eqref{eq:muest_lower}.  Alternatively, this bound can also be proven by the following method similar to the one in \cite{Kadets_Snobar}, which we reuse in Section~\ref{sec:extensions}.

Let $T$ be a simplex of maximal volume contained in $K^\circ$ with one vertex fixed at $o$. By the linear invariance property of $\mu(K,n)$, we may suppose that $T = \conv \{ o, e_1, \ldots, e_d \}$, where $\{e_j\}_{j=1}^d$ is the standard orthonormal basis of $\R^d$. Since the volume of $T$ is maximal, we derive that 
\[
K^\circ \subset B_\infty^d.
\]
Consequently, by taking polars, 
\[
    B_1^d
    \subset K.
\]
Therefore, for any $u \in \bd K$, 
\begin{equation}\label{eq:1norm_ineq}
    1
    = \|u\|_K
    \leq \| u \|_1
    = \sum_{j=1}^d |\la u, e_j\ra|.
\end{equation}
Thus, for any choice of $u_1, \ldots, u_n \in \bd K$,
\[
n \leq \sum_{i=1}^n \sum_{j=1}^d|\la u_i, e_j\ra| = \sum_{j=1}^d \Big(\sum_{i=1}^n|\la u_i, e_j\ra| \Big)
\]
which implies that
\[
\frac n d \leq \sum_{i=1}^n|\la u_i, e_j\ra|
\]
for some $j \in [d]$. Since $e_j \in K^\circ$, this shows that $\frac n d \leq  \max_{v \in K^\circ } \sum_{i=1}^n |\la u_i, v \ra|$. As the choice of $u_1, \ldots, u_n$ was arbitrary, we arrive at $\frac 1 d \leq \mu(K,n) $.
\end{remark}

\section{General Upper Bounds}
\label{sec:upper}

In this section, we verify the asymptotically sharp upper bound on the $n$-th polarization constant.
Doing so amounts to finding a set of unit vectors all of whose signed sums have sufficiently bounded norm.
Our sources of existence for these sets are the well-known Dvoretzky--Rogers lemma \cite{dvoretzky_rogers} and John ellipsoid theorem \cite{john1948extremum}.
While the former approach is used already in \cite{dvoretzky_rogers} to essentially upper bound $\mu(K)$,
the estimate we achieve is more precise and tailored to $\mu(K,n)$.

\begin{proof}[Proof of Theorem~\ref{thm:upperbd}]
First, we establish \eqref{eq:upperbd}.
By Lemma~\ref{lem:symm}, we may assume that $K$ is symmetric. The triangle inequality gives $\mu(K,n)\leq 1$. Thus for $d \leq 4$ or $n \leq 4$,
\[
    \mu(K,n)
    \leq 1
        \leq 2  \max \Big\{ \frac {1}{\sqrt{n}}, \frac{1}{\sqrt{d}} \Big\},
\]
so we assume $d,n \geq 5$ from now on.

By Lemma~\ref{lem:symm} and invariance under invertible linear transformations,
we may assume that $K$ is a symmetric convex body in John position; that is, $B_2^d$ is the maximal volume ellipsoid in $K$.
Consequently, for any $x \in \R^d$,
\begin{equation}\label{eq:Keucl}
    \|x\|_K
    \leq |x|.
\end{equation}
By the Dvoretzky--Rogers lemma \cite{dvoretzky_rogers},
there exist contact points $x_1, \ldots, x_d \in \bd K \cap S^{d-1}$ and an orthonormal basis $b_1, \ldots, b_d$ of $\R^d$ such that 
\begin{equation}\label{eq:DR}
    \la x_i, b_i \ra
    \geq \sqrt{\frac{d - i +1}{d}}
\end{equation}
for every $i \in [d]$. For each such index $i$, define
\[
    w_i
    := \frac 1 {\| b_i \|_K} \cdot b_i.
\]
Then $w_1, \ldots, w_d$ is an orthogonal system of vectors in $\bd K$. We build our choice of the $u_i$ from $w_1, \ldots, w_d$; in order to show that the $K$-norm of any signed sum is not too large, it suffices to estimate the Euclidean norm by \eqref{eq:Keucl}. Note that 
\[
K \subset \{ y \in \R^d: \la y, x_i \ra \leq 1 \text{ for all } i \in [d] \}. 
\]
Accordingly,
$\la w_i, x_i \ra \leq 1$ for each $i \in [d]$, so \eqref{eq:DR} yields
\begin{equation}\label{eq:winorm}
    |w_i|
    = \frac{\la w_i, x_i \ra}{\la b_i, x_i \ra}
    \leq \sqrt{\frac{d}{d - i +1}}.
\end{equation}

To demonstrate the method, we first prove for $n \geq d$ the weaker bound with the constant $2$ in \eqref{eq:upperbd} being replaced by $4$. We then refine the argument to achieve~\eqref{eq:upperbd}.

Now, let $k := \big\lceil \frac{d}{2} \big\rceil $ and set
$r := \big\lfloor \frac{n}{k} \big\rfloor$.
Then $ k \leq \frac 2 3 d$, and since we assumed $n \geq d$,
\[
    r + 1
    \leq \bigg\lfloor \frac{2n}{d} \bigg\rfloor + 1
    \leq \frac{3n}{d}.
\]
Let the vector set  $u_1, \ldots, u_m$ consist of $r+1$ 
copies of $w_1, \ldots, w_k$, where $m = k (r+1) \geq n$. Note that \eqref{eq:winorm} implies that $|w_i| < \sqrt{2}$ for all $i \in [k]$. Thus, by the orthogonality of $\{w_i\}_{i=1}^k$ and the above bounds on $k$ and $r+1$, 
it follows for any $\eps \in \{\pm1\}^m$ that
\[
    \Big| \sum_{i=1}^m \varepsilon_i u_i \Big|^2
    \leq \sum_{i=1}^k |(r+1) w_i|^2
    < 2 k (r+1)^2
    \leq 12 \frac{n^2}{d}.
\]
Therefore, \eqref{eq:Keucl} along with
the monotonicity of $n \cdot \mu(K,n)$ in $n$ shows that
\[
    n \cdot \mu(K,n)
    \leq m \cdot \mu(K,m)
    \leq \max_{\varepsilon \in \{\pm1\}^m} \Big\| \sum_{i=1}^m \varepsilon_i u_i \Big\|_K
    \leq \max_{\varepsilon \in \{\pm1\}^m} \Big| \sum_{i=1}^m \varepsilon_i u_i \Big|
    < 4 \frac{n}{\sqrt{d}}.
\]

To obtain the stronger upper bound \eqref{eq:upperbd}, we handle the number of copies of $w_1, \ldots, w_k$ more delicately and use more precise estimates.
Set $\ell \in \{0, \ldots, k-1\}$ so that $n = r \cdot k + \ell$.
We choose $u_1, \ldots, u_n$ to have $r+1$ copies of $w_1, \ldots, w_\ell$ and $r$ copies of $w_{\ell+1}, \ldots, w_k$.
By the orthogonality of $\{w_i\}_{i=1}^k$, \eqref{eq:winorm}, and the inequality $\frac{d}{d-i+1} < 2$ that holds for every $i \in [k]$,
it follows for any $\eps \in \{\pm1\}^n$ that
\begin{align}
\begin{split}
\label{eq:epsu1}
    \Big| \sum_{i=1}^n \varepsilon_i u_i \Big|^2
    & \leq \sum_{i=1}^\ell | (r+1) w_i |^2 + \sum_{i=\ell+1}^k | r w_i |^2
    \\
    & \leq \sum_{i=1}^\ell (r+1)^2 \Big( 2 + \frac{d}{d-i+1} - 2 \Big) + \sum_{i=\ell+1}^k r^2 \Big( 2 + \frac{d}{d-i+1} - 2 \Big)
    \\
    & \leq 2 \big( \ell (r+1)^2 + (k-\ell) r^2 \big) + \sum_{i=1}^k \Big( \frac{d}{d-i+1} - 2 \Big).
\end{split}
\end{align}
We estimate the first term above via
\begin{equation}\label{eq:epsu2}
    \ell (r+1)^2 + (k-\ell) r^2
    = k r^2 + 2 \ell r + \ell
    = \frac {n^2 } {k} - \frac{\ell^2}{k} + \ell
    \leq \frac {n^2 } {k} + \frac k 4
    \leq \frac {2 n^2 } {d} + \frac d 8,
\end{equation}
where the penultimate inequality is implied by the fact that the function $ f(\ell) = - \frac{\ell^2}{k} + \ell$ attains its maximum at $\ell = \frac k 2$, and the final inequality by the fact that $g(k) = \frac{n^2}{k} + \frac{k}{4}$ is decreasing on the interval $(0,2n)$. Hence, $k \geq \frac d 2$ shows that $g(k) \leq g(\frac d 2)$.

Next, note that 
\begin{equation}\label{eq:epsu3}
    \sum_{i=1}^k \Big( \frac{d}{d-i+1} - 2 \Big)
    \leq \int_{d-k}^d \frac{d}{t} \, \dd t - 2 k
    = d \ln \Big( \frac{d}{d-k} \Big) - 2k.
\end{equation}
If $d$ is even, then $k = \frac{d}{2}$ and accordingly
\[
    d \ln \Big( \frac{d}{d-k} \Big) - 2 k
    = d (\ln(2) - 1)
    < - \frac{d}{4}.
\]
If instead $d$ is odd, so that $k = \frac{d+1}{2}$, then
\begin{align*}
    d \ln \Big( \frac{d}{d-k} \Big) - 2 k
    & = d \ln \Big( \frac{2d}{d-1} \Big) - d-1
    = d \ln \Big( 2 + \frac{2}{d-1} \Big) - d - 1
    \\
    & \leq d \Big( \ln(2) + \frac{1}{d-1} \Big) - d - 1
    = d (\ln(2)-1) + \frac{1}{d-1},
\end{align*}
where we used the standard inequality $\ln (1 + x)  \leq x $ that holds for $x \geq 0$.
With $d \geq 5$,
\[
    d \ln \Big( \frac{d}{d-k} \Big) - 2 k
    \leq d (\ln(2)-1) + \frac{1}{d-1}
    \leq d (\ln(2)-1) + \frac{d}{20}
    < - \frac{d}{4}.
\]

In either case, combining the above inequalities with \eqref{eq:epsu1}, \eqref{eq:epsu2}, and \eqref{eq:epsu3} yields that
\[
    \Big| \sum_{i=1}^n \varepsilon_i u_i \Big|^2
    < 4 \frac{n^2}{d}.
\]
Therefore, if $n \geq d$, then by \eqref{eq:Keucl},
\[
    \mu(K,n)
    \leq \max_{\varepsilon \in \{\pm1\}^n} \frac 1 n \Big\| \sum_{i=1}^n \varepsilon_i u_i \Big\|_K
    \leq \max_{\varepsilon \in \{\pm1\}^n} \frac 1 n \Big| \sum_{i=1}^n \varepsilon_i u_i \Big|
    < \frac{2}{\sqrt{d}}.
\]

Finally, assume that $5 \leq n <d$. Take any $n$-dimensional linear subspace $E \subset \R^d$ and apply the above argument to $K \cap E$ within $E$.
This yields some $u_1, \ldots, u_n \in \bd K \cap E$ with
\[
   \mu(K,n)
    \leq \max_{\varepsilon \in \{\pm1\}^n} \frac 1 n \Big\| \sum_{i=1}^n \varepsilon_i u_i \Big\|_K
    < \frac{2}{\sqrt{n}},
\]
thus completing the proof of \eqref{eq:upperbd}.

Next, we turn to \eqref{eq:upperbd2}. We may assume, as before,  that $K \in \K_o^d$ is a symmetric convex body in John position. Accordingly, there exist some points $x_1, \ldots, x_N \in \bd K \cap S^{d-1}$ with $N \leq \frac{d (d+1)}{2}$ and weights $\lambda_1, \ldots, \lambda_N > 0$ with $\sum_{i=1}^N \lambda_i = d$ such that
\[
    \sum_{i=1}^N \lambda_i \la v, x_i \ra x_i
    = v
\]
for any $v \in \R^d$.
Now, for any given $n \geq 1$,
take $m_i := \lceil \frac{n}{d} \cdot \lambda_i \rceil$ copies of $x_i$ for all $i \in [N]$ to obtain some vectors $u_1, \ldots, u_m \in \bd K$
with $m = \sum_{i=1}^N \lceil \frac{n}{d} \cdot \lambda_i \rceil \geq \frac{n}{d} \cdot \sum_{i=1}^N \lambda_i = n$.
Since $m_i < \frac{n}{d}\cdot\lambda_i+1$ and
$|\la x_i,v\ra|\leq 1$ for $v\in K^\circ\subset B_2^d$,
the Cauchy--Schwarz inequality yields
\begin{align*}
    \sum_{i=1}^m |\la u_i, v \ra|
    & = \sum_{i=1}^N m_i |\la x_i, v \ra|
    < N + \sum_{i=1}^N \frac{n}{d} \cdot \lambda_i|\la x_i, v \ra|
    \\
    & \leq \frac{d (d+1)}{2} + \frac{n}{d} \cdot \sqrt{ \sum_{i=1}^N \lambda_i } \cdot \sqrt{ \sum_{i=1}^N \lambda_i \la x_i, v \ra^2 }
    \\
    & = \frac{d (d+1)}{2} + \frac{n}{d} \cdot \sqrt{d} \cdot |v|
    \leq \frac{d (d+1)}{2} + \frac{n}{\sqrt{d}}.
\end{align*}
This shows
\[
    n \cdot \mu(K,n)
    \leq m \cdot \mu(K,m)
    < \frac{d (d+1)}{2} + \frac{n}{\sqrt{d}}
\]
and, consequently,
\[
    \mu(K,n)
    < \frac{1}{\sqrt{d}} + \frac{d (d+1)}{2n}.
\]
This establishes \eqref{eq:upperbd2}, and completes the proof of Theorem~\ref{thm:upperbd}.
\end{proof}

\begin{remark}
A weaker, asymptotically matching upper bound can be derived from a result of Giannopoulos~\cite{giannopoulos1996proportional}, which states that for any symmetric $K \in \K^d_o$ and any $\varepsilon \in (0,1)$,
there exist vectors $x_1, \ldots, x_m \in \R^d$, $m \geq (1 - \eps) d$, such that for any reals $t_1, \ldots, t_m \in \R$, 
\begin{equation*}
    \max_{i \leq m} |t_i|
    \leq \Big \| \sum_{i=1}^m t_i x_i \Big \|_K \leq  a \eps^{-b} \Big( \sum_{i=1}^m t_i^2\Big)^{1/2},
\end{equation*}
where $a,b$ are positive constants. For the present application, one has to set $\eps = \frac 12 $ and build the vector set from rescaled copies of $x_1, \ldots, x_m$.
\end{remark}

\section{Limit Transition}\label{sec:Limit}

Next, we prove that for any $K \in \K_o^d$,
the $n$-th large signed sum and polarization constants of $K$ converge to its Macphail and $1$-absolutely summing constants as $n \to \infty$.
By the proof of Theorem~\ref{thm:muM}, it suffices to consider the polarization variants.
Alternatively, the proof given below can be adjusted to work for the signed sum setting as well.
The main intuition behind the proof is that with growing numbers of unit vectors and combining copies of the same vector into a single longer one, we can approximate any set of vectors arbitrarily well up to scaling.

\begin{proof}[Proof of Theorem~\ref{thm:pollim}]
\newcommand\multvar{k}
We first show that the limit exists, following the proof of \cite[Proposition~14.6.1]{borodachov2019discrete}.
Fix $K \in \K_o^d$, and for each $i \geq 1$, let $U_i$ be an $i$-element vector family contained in $\bd K$ with
\[
    \mu(K,i)
    = \max_{v \in K^\circ} \frac 1 i \sum_{u \in U_i} |\la u, v \ra |.
\]
Let now $n,m \geq 1$ be arbitrary, and set $U = U_n \cup U_m$ (understood as the union of families rather than sets).
Then
\begin{align}
\begin{split}\label{eq:mu_subadd}
    (n+m) \cdot \mu(K,n+m)
    & \leq  \max_{v \in K^\circ} \sum_{u \in U} |\la u, v \ra |
    \leq  \max_{v \in K^\circ} \sum_{u \in U_n} |\la u, v \ra | + \max_{w \in K^\circ} \sum_{u \in U_m} |\la u, w \ra | \\ 
    & = n \cdot \mu(K,n) + m \cdot \mu(K,m),
\end{split}
\end{align}
meaning the sequence $\{ n \cdot \mu(K,n): n = 1, 2, \ldots \} $ is subadditive. Consequently, by Fekete's lemma for subadditive sequences \cite{fekete1923verteilung},
\[
\lim_{n \to \infty} \mu(K,n) = \inf_{n \geq 1} \mu(K,n).
\]
Therefore,
\begin{align*}
    \lim_{n \to \infty} \mu(K,n)
    & = \inf_{ U \subset \bd K \text{ finite family}} \ \max_{v \in K^\circ} \frac{\sum_{u \in U} |\la u, v \ra|}{ |U| } \\
    & = \inf_{ U \subset \bd K \text{ finite family}} \ \max_{v \in K^\circ} \frac{\sum_{u \in U} |\la u, v \ra|}{ \sum{_{u \in U}} \| u\|_K }.
\end{align*}

Next, we show that the latter quantity remains unchanged if we extend the domain of the infimum from finite vector families $U \subset \bd K$ to finite vector families $X \subset \R^d \setminus \{ o \}$. According to the remark following Definition~\ref{def:macphail}, one can then restrict the domain to finite vector sets $X \subset \R^d \setminus \{ o \}$, and the result follows. 

We use the idea of the proof of \cite[Theorem~2.2]{DL23}. For a family $X \subset \R^d \setminus \{o\}$, let
\[
    m(X)
    := \max_{v \in K^\circ} \frac{\sum_{x \in X} |\la x, v \ra|}{ \sum{_{x \in X}} \| x\|_K }.
\]
Clearly, 
\[
    \inf_{ U \subset \bd K \text{ finite family}} m(U)
    \geq \inf_{ X \subset \R^d \setminus \{ o \} \text{ finite family}} m(X),
\]
so we only have to prove the reverse inequality. Let $\eps >0$ be arbitrary, and let $X_\eps$ be a finite vector family for which $m(X_\eps)$ differs from the right-hand infimum by less than $\varepsilon$.
After an arbitrarily small perturbation of $X_\eps$ if necessary, we can assume for every $x \in X_\eps$ that $\| x \|_K$ is rational.
Let $\multvar$ be the smallest common multiple of the denominators of these norms. Replacing each $x \in X_\eps$ with $\multvar \| x \|_K$ copies of $\frac{1}{\multvar \| x \|_K} \cdot x $ leads to a vector family $X_\eps'$ with $m(X_\eps) = m(X_\eps')$.
Moreover, $X_\eps'$ is contained in $\frac 1 \multvar \bd K$.
Scaling by a factor of $\multvar$ leads to a family in $\bd K$ with still the same $m$-value, and letting $\eps \to 0$ finishes the proof. 
\end{proof}

\section{Extremal Configurations in Special Cases} \label{sec:special}

In this section, we prove our results on extremal configurations in low dimensions, for small values of $n$, and for certain special families of vectors, thereby resolving particular cases of Problems~\ref{problem1} and~\ref{problem2}. In particular, we prove the nearly sharp planar upper bounds stated in Corollary~\ref{cor:planar_upper}, as well as the sharp lower bound for the convex bodies in Proposition~\ref{prop:etf}. We begin with the latter since we use it to establish the sharpness of the former.

\begin{proof}[Proof of Proposition~\ref{prop:etf}]
As a first step, we determine the parameters $\alpha$ and $\beta$ that appear in the definition of the ETF and verify in passing that the set range specified in \eqref{eq:Kpolcond} is indeed non-empty. By taking traces in \eqref{eq:TF}, we immediately derive that 
\[
\alpha = \frac N d.
\]
Next, choose $j \in [N]$ and let $\varepsilon_i$ be the sign of $\la w_j, w_i \ra$ for each $i \in [N]$, or arbitrary when the inner product is 0.
Then, by \eqref{eq:TF1}, 
\begin{equation}\label{eq:wj_cont}
    \frac{N}{d} w_j
    = \sum_{i=1}^N \langle w_j, w_i \rangle w_i
    = w_j + \beta \sum_{i \neq j} \varepsilon_i w_i.
\end{equation}
Taking the inner product of both sides with $w_j$ and rearranging leads to
\[
   \beta
   = \sqrt{\frac{\etfind-d}{d (\etfind-1)}}.
\]
If now $N = d$, then $\gamma_{d,d} = 1$ yields $w_j \in \frac{1}{\gamma_{d,N}} \sum_{i=1}^N [-w_i,w_i]$.
Otherwise, by \eqref{eq:wj_cont},
\[
    w_j
    = \frac{\beta}{\beta + \frac{N}{d} - 1} \sum_{i=1}^N \varepsilon_i w_i
    = \frac{1}{1 + \frac{N-d}{\beta d}} \sum_{i=1}^N \varepsilon_i w_i
    = \frac{1}{\gamma_{d,N}} \sum_{i=1}^N \varepsilon_i w_i
    \in \frac{1}{\gamma_{d,N}} \sum_{i=1}^N [-w_i,w_i],
\]
so we conclude from the symmetry of the right-hand zonotope that the range specified in \eqref{eq:Kpolcond} is indeed non-empty.

For the rest of the proof, by using the invariance of $\mu$ under linear transformations and the symmetry property provided by Lemma~\ref{lem:symm}, we may assume that $K$ is symmetric with 
\begin{equation*}
    \conv \{\pm \etfvar_1,\ldots, \pm \etfvar_\etfind\} \subset K^\circ
    \subset \frac 1 {\gamma_{d,\etfind}} \sum_{i=1}^{\etfind} [-\etfvar_i,\etfvar_i].
\end{equation*}
For any $u_1, \ldots, u_n \in \bd K$, by \eqref{eq:polarsupport},
\begin{align*}
\begin{split}
    n
    & = \sum_{i=1}^n \| u_i \|_K
    = \sum_{i=1}^n h_{K^\circ}(u_i)
    \leq \frac 1 {\gamma_{d,\etfind}} \sum_{i=1}^n h_{\sum_j[-\etfvar_j,\etfvar_j]}(u_i) \\
    & = \frac 1 {\gamma_{d,\etfind}} \sum_{i=1}^n \sum_{j=1}^N|\la u_i, w_j \ra |
    = \frac 1 {\gamma_{d,\etfind}} \sum_{j=1}^N  \sum_{i=1}^n  |\la u_i, w_j \ra |
    \leq  \frac 1 {\gamma_{d,\etfind}} N \max_{v \in K^\circ} \sum_{i=1}^n  |\la u_i, v \ra |,
\end{split}
\end{align*}
where in the last inequality we used that $w_j \in K^\circ$ for every $j \in [N]$. Consequently, 
\[
    \mu(K,n)
    = \min_{u_1, \ldots, u_n \in \bd K} \ \max_{v \in K^\circ } \frac 1 n \sum_{i=1}^n |\la u_i, v \ra |
    \geq \frac 1 n \cdot n \cdot \frac{\gamma_{d, N}}{N}
    = \frac{\gamma_{d, N}}{N}.
\]

For the claimed equality case, we assume that $K^\circ = \conv \{ \pm w_1, \ldots, \pm w_N \}$ and $n = k N$ for some integer $k \geq 1$. Since $|w_i| = 1$ for every $i \in [N]$, the assumption on $K$ yields $w_i \in \bd K$.
We can therefore choose $u_1, \ldots, u_n$ to consist of $k$ copies of $w_1, \ldots w_N$.
Then, for any $j \in [\etfind]$,
\[
    \sum_{i=1}^n |\la u_i, \etfvar_j \ra|
    = k \sum_{i=1}^N |\la w_i, w_j \ra|
    = k (1 + (N-1) \beta)
    = k \Bigg( 1 + \sqrt{ \frac{(\etfind-d) (\etfind-1)}{d} } \Bigg)
    = k \gamma_{d,N},
\]
hence, by Lemma~\ref{lem:convhull},
\[
    \mu(K,n) 
    = \mu(K,k \etfind)
    \leq \max_{v \in K^\circ} \frac{1}{n} \sum_{i=1}^n |\la u_i, v \ra|
    = \max_{j \in [\etfind]} \frac{1}{n} \sum_{i=1}^n |\la u_i, \etfvar_j \ra|
    = \frac{\gamma_{d,N}}{N}. \qedhere
\]
\end{proof}

Having established all the ingredients needed for Theorem~\ref{thm:mu_upper}, we now combine them into a short proof.

\begin{proof}[Proof of Theorem~\ref{thm:mu_upper}]
By  Lemma~\ref{lem:symm}, we may assume that $K$ is symmetric.
The upper estimate on $\mu(K)$ follows from Propositions~\ref{prop:gordon}~and~\ref{prop:proj_strong}.
If $K$ satisfies the criterion \eqref{eq:mumax_characterization} for a maximal ETF, then Theorem~\ref{thm:pollim} and Proposition~\ref{prop:etf} verify the matching lower bound on $\mu(K)$.
If $K$ does not satisfy \eqref{eq:mumax_characterization} for a maximal ETF, then $\lambda(K) < d \gamma_d$ from Proposition~\ref{prop:proj_strong} together with Propositon~\ref{prop:gordon} shows that equality in \eqref{eq:mu_uppergamma} cannot occur.
\end{proof}

We move on to the sharp upper bounds in Propositions~\ref{prop:n=2}~and~\ref{prop:planar_n=3}.
Our approaches are based on inscribing affine-regular polygons in convex curves.

\begin{proof}[Proof of Proposition~\ref{prop:n=2}]
By Lemma~\ref{lem:symm}, we may assume that $K$ is symmetric.
If $d=1$, then clearly $\mu(K,2) = 1$.
Otherwise, let $E \subset \R^d$ be any $2$-dimensional linear subspace.
By a result from \cite{gruenbaum1959}, there exists a symmetric affine-regular octagon $O$ inscribed in $K \cap E$.
Choose vertices $u_1,u_2$ of $O$ such that none of $\pm u_1, \pm u_2$ are adjacent vertices of $O$.
Then $\frac{1}{\sqrt{2}} (\varepsilon_1 u_1 + \varepsilon_2 u_2)$ for $\varepsilon \in \{\pm 1\}^2$ is always a vertex of $O$ by affine-regularity.
Since $u_1,u_2 \in \bd K$, we obtain from
\eqref{eq:euclidean_small_n} that
\[
    \mu(K,2)
    \leq \max_{\varepsilon \in \{\pm 1\}^2} \frac 1 2 \| \varepsilon_1 u_1 + \varepsilon_2 u_2 \|_K
    = \frac 1 {\sqrt{2}}
    = \mu(B_2^d,2). \qedhere
\]
\end{proof}

\begin{proof}[Proof of Proposition~\ref{prop:planar_n=3}]
Theorem~\ref{thm:sharplower} yields $\mu(K,3) \geq \frac 2 3$.
To show $\mu(K,3) \leq \frac 2 3$, we assume that $K$ is symmetric by Lemma~\ref{lem:symm} and take a symmetric affine-regular hexagon $H$ inscribed in $K$ (see, e.g., \cite{gruenbaum1959}).
Let $\pm u_1, \pm u_2, \pm u_3$ be the vertices of $H$ and $\pm v_1, \pm v_2, \pm v_3$ the vertices of $H^\circ$ such that $|\la u_i, v_j \ra| = 1$ for $i,j \in [3]$, $i \neq j$.
The affine-regularity of $H$ yields $\la u_i, v_i \ra = 0$ for all $i \in [3]$.
Since $u_1, u_2, u_3 \in \bd K$, we obtain from $K^\circ \subset H^\circ$ and Lemma~\ref{lem:convhull} that
\[
    \mu(K,3)
    \leq \max_{w \in K^\circ} \frac 1 3 \sum_{i=1}^3 |\la u_i, w \ra|
    \leq \max_{v \in H^\circ} \frac 1 3 \sum_{i=1}^3 |\la u_i, v \ra|
    = \max_{j \in \{1,2,3\}} \frac 1 3 \sum_{i=1}^3 |\la u_i, v_j \ra|
    = \frac 2 3 .\qedhere
\]
\end{proof}

With all the above results verified,
we obtain Corollary~\ref{cor:planar_upper} from \eqref{eq:mu_subadd}.

\begin{proof}[Proof of Corollary~\ref{cor:planar_upper}]
We begin with the upper bound when $n \equiv 0$ (mod $3$).
Induction over $k \geq 1$, with Proposition~\ref{prop:planar_n=3} covering the base case $k=1$, shows with \eqref{eq:mu_subadd} that
\[
    \mu(K,3(k+1))
    \leq \frac{3k}{3 (k+1)} \cdot \mu(K,3k) + \frac{3}{3(k+1)} \cdot \mu(K,3)
    \leq \frac{k}{k+1} \cdot \frac{2}{3} + \frac{1}{k+1} \cdot \frac{2}{3}
    = \frac{2}{3}.
\]
The upper bound for $n = 2$ is given by Proposition~\ref{prop:n=2}.
Consequently, if $n \geq 5$ and $n \equiv 2$ (mod $3$), then \eqref{eq:mu_subadd} leads to
\[
    \mu(K,n)
    \leq \frac{n-2}{n} \cdot \mu(K,n-2) + \frac{2}{n} \cdot \mu(K,2)
    \leq \frac{n-2}{n} \cdot \frac{2}{3} + \frac{2}{n} \cdot \frac{1}{\sqrt{2}}
    = \frac{2}{3} + \frac{1}{n} \Big( \sqrt{2} - \frac{4}{3} \Big).
\]
For the final set of inequalities, Proposition~\ref{prop:n=2} with \eqref{eq:mu_subadd} gives
\[
    \mu(K,4)
    \leq \frac{2}{4} \cdot \mu(K,2) + \frac{2}{4} \cdot \mu(K,2)
    \leq \frac{1}{\sqrt{2}}.
\]
For $n \geq 7$ with $n \equiv 1$ (mod $3$), using \eqref{eq:mu_subadd} again thus leads to
\[
    \mu(K,n)
    \leq \frac{n-4}{n} \cdot \mu(K,n-4) + \frac{4}{n} \cdot \mu(K,4)
    \leq \frac{n-4}{n} \cdot \frac{2}{3} + \frac{4}{n} \cdot \frac{1}{\sqrt{2}}
    \leq \frac{2}{3} + \frac{2}{n} \Big( \sqrt{2} - \frac{4}{3} \Big).
\]
The sharpness of the inequality for $n \leq 3$ follows from Propositions~\ref{prop:n=2}~and~\ref{prop:planar_n=3},
and for $n \equiv 0$ (mod $3$) from Proposition~\ref{prop:etf}.
The required ETF with three vectors comes from the vertices of a regular hexagon inscribed in $S^1$.
\end{proof}

To close this section, let us point out that the upper bound in Corollary~\ref{cor:planar_upper} is nearly sharp for all $n \geq 1$,
as Proposition~\ref{prop:etf} yields that for a symmetric regular hexagon $H \in K_o^2$,  the bound $\mu(H,n) \geq \frac{2}{3}$ holds for all $n \geq 1$.
In fact, it is not hard to prove that $\mu(H,n) = \frac{2}{3}$ for all $n \geq 2$.

\section{The Euclidean Case}\label{sec:Euclidean}

In this section, we prove the results regarding the Euclidean ball.
The sharp lower and upper bounds in Proposition~\ref{prop:polball} follow essentially from the results and methods in \cite[Theorem~1~\&~Theorem~2]{AN2020}
(see also \cite[Theorem~2]{Gordon69} together with Theorem~\ref{thm:pollim} for the lower bound).
To be self-contained, we describe a short proof of the upper bound based on tight frames.
We present an alternative approach for the lower bound by establishing a connection to the first intrinsic volume and taking advantage of its Minkowski linearity.
We further sharpen this approach in the planar case to obtain a direct proof of Proposition~\ref{prop:planar}.

\begin{proof}[Proof of the upper bound in Proposition~\ref{prop:polball}]
For any $n \geq d$, let $u_1, \ldots, u_n \in S^{d-1}$ form a tight frame.
The existence of such vector sets for all $n \geq d$ is proven in \cite{benedetto2003finite}.
Considering traces for the exact constant $\alpha$ of the tight frame, we have
\[
    \sum_{i=1}^n \la x, u_i \ra u_i
    = \frac{n}{d} \cdot x
\]
for any $x \in \R^d$. 
Then, for any $v \in S^{d-1}$, 
\[
    \sum_{i=1}^n |\la u_i, v \ra|^2
    = \frac n d.
\]
By the Cauchy--Schwarz inequality, this leads to
\[
    \sum_{i=1}^n |\la u_i, v \ra|
    \leq \sqrt{n} \cdot \sqrt{\sum_{i=1}^n |\la u_i, v \ra|^2} = \frac {n}{\sqrt{d}},
\]
which implies $\mu(B_2^d,n) \leq \frac 1 {\sqrt{d}}$.
\end{proof}

To derive the  lower bound in~\eqref{eq:ballpol}, we invoke some classical geometric inequalities.
We write $R(K)$ for the circumradius of a non-empty compact convex set $K \subset \R^d$, i.e., the smallest possible radius of a ball containing $K$. Furthermore, let $V_1(K)$ denote the first intrinsic volume of $K$ (see, e.g., \cite{2020BöröczkyHug,gardner1995geometric, schneider2013convex}) -- then $ \frac {2 \, \kappa_{d-1}}{d \, \kappa_d}V_1(K)$ is the mean width of $K$:
\begin{equation}\label{eq:meanwidth}
    V_1(K)
    = \frac{d \, \kappa_d} {2 \, \kappa_{d-1}} \int_{S^{d-1}} h_K(u) + h_K(-u) \, \dd \sigma(u),
\end{equation}
where $\sigma(.)$ denotes the rotation-invariant probability measure on $S^{d-1}$.
The circumradius and the first intrinsic volume are linked by the inequalities
\begin{equation}\label{eq:R_V1}
    \frac{d \, \kappa_d}{\kappa_{d-1}} R(K)
    \geq V_1(K)
    \geq 2 R(K),
\end{equation}
where the left-hand inequality follows from
the fact that the width of $K$ in any given direction is at most $2 R$,
and the right-hand inequality is shown in \cite[Theorem~1.4]{2020BöröczkyHug}.
The left-hand inequality is sharp only for balls,
and the right-hand inequality becomes an equality precisely for segments on the class of non-empty compact convex sets.

\begin{proof}[Proof of the lower bound in Proposition~\ref{prop:polball}]
For $u_1, \ldots, u_n \in S^{d-1}$, define the zonotope $Q = Q(u_1, \ldots, u_n):= \sum_{i=1}^n [-u_i,u_i]$.
It is clear that
\begin{equation}\label{eq:zonotope_LSS}
    R(Q)
    = \max_{v \in Q} |v|
    = \max_{\varepsilon \in \{\pm1\}^n} \Big| \sum_{i=1}^n \varepsilon_i u_i \Big|,
\end{equation}
while the Minkowski linearity of the mean width, which is apparent from \eqref{eq:meanwidth}, directly implies that
\[
    V_1(Q)
    = \sum_{i=1}^n V_1([-u_i, u_i]).
\]
Therefore, \eqref{eq:R_V1} yields
\begin{align}
\begin{split}\label{eq:zonotope_radius}
    R(Q)
    & \geq \frac{\kappa_{d-1}}{d \, \kappa_d} V_1(Q)
    = \frac{\kappa_{d-1}}{d \, \kappa_d} \sum_{i=1}^n V_1([-u_i,u_i]) \\
    & = \frac{2 \, \kappa_{d-1}}{d \, \kappa_d} \sum_{i=1}^n R([-u_i,u_i])
    = \frac{2 \, \kappa_{d-1}}{\sqrt{d} \, \kappa_d} \cdot \frac{n}{\sqrt{d}}.
\end{split}
\end{align}
Accordingly, 
\[
\mu(B_2^d,n) = \min_{u_1, \ldots, u_n \in S^{d-1}} \frac{R\bigl(Q(u_1, \ldots, u_n)\bigr)}{n} \geq  \frac{2 \, \kappa_{d-1}}{\sqrt{d} \, \kappa_d} \cdot \frac{1}{\sqrt{d}}.
\qedhere
\]
\end{proof}
Choosing the $u_i$ such that $Q$ is nearly a ball shows that the inequality above is sharp as $n \to \infty$, i.e.,
\[
    \lim_{n \to \infty} \sqrt{d} \cdot \mu(B_2^d,n) 
    = \frac{2 \, \kappa_{d-1}}{\sqrt{d} \, \kappa_d}
    = \sqrt{\frac{2}{\pi}} + O \Big(\frac{1}{d} \Big).
\]

Moreover, allowing arbitrary norms of the $u_i$ and using $R([-u_i,u_i]) = |u_i|$ in \eqref{eq:zonotope_radius} gives a direct proof of Corollary~\ref{cor:eucl_general_vectors}.
The equality case follows since $Q$ is a (non-degenerate) ball if and only if $d=1$.

Our new proof of Proposition~\ref{prop:planar} is based on the above approach:
we use an improved version of the left-hand inequality in \eqref{eq:R_V1} specific to polygons, which is a consequence of the following isoperimetric principle for the circumradius of polygons.
We write $\per(K) = 2V_1(K)$ for the perimeter of a convex body $K \subset \R^2$.

\begin{proposition}\label{prop:dowker}
For any $n \geq m \geq 1$ and any convex $m$-gon $Q \subset \R^2$,
\[
	2 n \sin \bigg( \frac{\pi}{n} \bigg) R(Q)
	\geq \per(Q),
\]
with equality if and only if $Q$ is a regular $n$-gon or a singleton.
\end{proposition}
\noindent 
Above, (regular) $1$-gons and $2$-gons are understood as singletons and segments, respectively.

The proposition falls into the category of \emph{Dowker-type results}
(in reference to the fundamental paper \cite{dowker1944} by Dowker;
see \cite{eggleston1957,fejestoth1955,molnar1955} for the corresponding results involving the perimeter).
Despite the elementary nature of the problem,
a transparent proof of Proposition~\ref{prop:dowker} appears difficult to locate in the literature
(see \cite{basitlangi2024,fejestoths1973} for the inequality without the characterization of the equality case,
and \cite[Open~Problems~(3a)]{audethansenmessine2007} for the entire result but without proof).
For the sake of completeness, we provide a short proof.

\begin{proof}
Let $\nu := R(Q)$, and assume that $t + \nu B_2^2$ is the circumcircle of $Q$ with some $t  \in \R^2$.
All claims are clear if $Q$ is a singleton, so we may assume $\nu > 0$.
It is well-known that $t$ must be a non-extreme point of $Q$ in this case.
If $Q'$ denotes the convex $m$-gon whose vertices are the intersection points of the circle $t + \nu S^1$ with the rays emanating from $t$ through the vertices of $Q$, then $t \in Q \subset t + \nu B_2^2$ implies $Q \subset Q'$.
We may further enlarge $Q'$ to a convex $n$-gon $Q^*$ all of whose vertices lie on the circle $t + \nu S^1$.
The strict monotonicity of the perimeter among convex sets with respect to inclusion
shows $\per(Q) \leq \per(Q^*)$,
with equality if and only if $Q$ is an $n$-gon with all of its vertices on $t + \nu S^1$.
Since $R(Q^*) = R(Q)$, it suffices to prove the proposition for $Q^*$.

The assumption $\nu > 0$ implies by $R(Q^*) = \nu$ that $n \geq 2$.
Let $v_1, \ldots, v_n$ be the distinct vertices of $Q^*$ in circular order and write $v_{n+1} := v_1$.
For $i \in [n]$, let $\alpha_i \in [0,\pi]$ be the angle enclosed by the vectors $v_i - t$ and $v_{i+1} - t$. Since the vertices are on the circle of radius $\nu$,  the distance between $v_i$ and $v_{i+1}$ 
equals $2 \nu \sin ( \frac{\alpha_i}{2} )$.
Note that $t \in Q^*$ implies $\sum_{i=1}^n \alpha_i = 2 \pi$.
Hence, the concavity of sine on the interval $[0,\pi]$ and Jensen's inequality yield
\[
	\per(Q^*)
	= \sum_{i=1}^n 2 \nu \sin \bigg( \frac{\alpha_i}{2} \bigg)
	= 2 n \nu \sum_{i=1}^n \frac{\sin( \frac{\alpha_i}{2} )}{n}
	\leq 2 n \nu \sin \Bigg( \sum_{i=1}^n \frac{\alpha_i}{2n} \Bigg)
	= 2 n \sin \bigg( \frac{\pi}{n} \bigg) R(Q^*).
\]
Equality holds if and only if every $\alpha_i$ equals $\frac{2 \pi}{n}$,
equivalently, $Q^*$ is a regular $n$-gon as claimed.
\end{proof}

\begin{proof}[Proof of Proposition~\ref{prop:planar}]
Let $u_1, \ldots, u_n \in S^1$
and define the symmetric polygon $Q := \sum_{i=1}^n [-u_i,u_i]$. 
Since the perimeter is Minkowski additive, 
\begin{equation}\label{eq:perimeter}
	\per(Q)
	= \sum_{i=1}^n \per([-u_i,u_i])
    = \sum_{i=1}^n 4 |u_i|
    = 4 n.
\end{equation}
Furthermore, any edge of $Q$ is parallel to one of the segments $[-u_i,u_i]$.
Therefore, $Q$ has at most $2n$ edges (fewer if some $u_i$ coincide up to sign). Thus, by \eqref{eq:zonotope_LSS}, Proposition~\ref{prop:dowker}, and \eqref{eq:perimeter},
\begin{equation}\label{eq:planar}
	\max_{\varepsilon \in \{\pm1\}^n} \Big| \sum_{i=1}^n \varepsilon_i u_i \Big|
	= R(Q)
	\geq \frac{\per(Q)}{4 n \sin ( \frac{\pi}{2n} )}
	= \frac{1}{\sin ( \frac{\pi}{2n} )}.
\end{equation}
Taking the minimum over the $u_i$ 
yields $\mu(B_2^2,n) \geq \frac{1}{n  \sin ( \frac{\pi}{2n} )}$.

Equality holds in \eqref{eq:planar} if and only if $Q$ is a regular $2n$-gon. 
Note that in this case, the sides of $Q$ are precisely the translated copies of $[-u_i, u_i]$ for $i=1, \ldots, n$, each with multiplicity 2. Hence,  $Q$ is a regular $2n$-gon if and only if 
$\{\pm u_1, \ldots, \pm u_n\}$ forms the vertex set of a regular $2n$-gon.
This shows that $\mu(B_2^2,n) = \frac{1}{n \sin ( \frac{\pi}{2n} )}$, as well as the characterization of the extremal case.
\end{proof}

As before for Corollary~\ref{cor:eucl_general_vectors},
the above proof easily extends to families of not necessarily unit vectors.
The only difference is the omission of the last equality in \eqref{eq:perimeter}.

\begin{corollary}
For any $u_1, \ldots, u_n \in \R^2$,
\[
    \max_{\varepsilon \in \{\pm1\}^n} \Big| \sum_{i=1}^n \varepsilon_i u_i \Big|
    \geq \frac{1}{n \sin( \frac{\pi}{2n})} \sum_{i=1}^n |u_i|,
\]
with equality if and only if all $u_i$ are zero or $\{\pm u_1, \ldots, \pm u_n\}$ forms the vertex set of a regular $2n$-gon. 
\end{corollary}

\section{Extensions}
\label{sec:extensions}

We conclude the article with two generalizations of some of the preceding results. First, we replace Minkowski norms by support functions of general compact convex sets.
For $K \in \K^d_o$, \eqref{eq:polarsupport} states that $\|.\|_K = h_{K^\circ}$ is such a support function, but of a convex set that contains the origin $o$ in its interior.
Therefore, the subsequent results can be understood as the extension of Theorem~\ref{thm:sharplower} for convex bodies that do not possess this extra property.

\begin{remark}\label{rem:sublin_reduction}
It is straightforward to reduce from support functions $h_K$ of non-empty compact convex sets $K \subset \R^d$ back to Minkowski norms if we want to show, for example, inequalities of the form
\begin{equation}\label{eq:general_sublin_ineq}
    \max_{\varepsilon \in \{\pm1\}^n} h_K \Big( \sum_{i=1}^n \varepsilon_i u_i \Big)
    \geq c \sum_{i=1}^n |h_K(u_i)|
\end{equation}
for vectors $u_1, \ldots, u_n \in \R^d$ and constants $c = c(d,n) \geq 0$.
Define the function
\[
    g \colon \R^d \to \R,
    g(x) = \max \{ h_K(x), h_K(-x) \}
    = \max \{ h_K(x), h_{-K}(x)\}
    = h_{\conv \{K \cup (-K)\}}(x),
\]
which is a norm on the linear hull of $K$ by \eqref{eq:polarsupport}. Further note for $x in \R^d$ that
\[
    0
    = h_K(o)
    = h_K(x-x)
    \leq h_K(x) + h_K(-x)
\]
from the properties of $h_K$ gives
\begin{equation}\label{eq:max_sublin}
    g(x)
    = \max \{ h_K(x), h_K(-x) \}
    \geq \max \{ h_K(x), -h_K(x) \}
    = |h_K(x)|.
\end{equation}
The orthogonal projection $v$ of some $u \in \R^d$ onto the linear hull of $K$ satisfies $h_K(v) = h_K(u)$ and $g(v) = g(u)$.
Thus, with $v_i$ being the orthogonal projection of $u_i$ for $i \in [n]$,
\[
    \max_{\varepsilon \in \{\pm1\}^n} h_K \Big( \sum_{i=1}^n \varepsilon_i u_i \Big)
    = \max_{\varepsilon \in \{\pm1\}^n} h_K \Big( \sum_{i=1}^n \varepsilon_i v_i \Big)
    = \max_{\varepsilon \in \{\pm1\}^n} g \Big( \sum_{i=1}^n \varepsilon_i v_i \Big).
\]
Together with \eqref{eq:max_sublin}, we see that \eqref{eq:general_sublin_ineq} for $g$ and the $v_i$ implies \eqref{eq:general_sublin_ineq} for $h_K$ and the $u_i$ with the same constant $c$.
We may further assume that the subspace spanned by the $v_i$, which is a subspace of the linear hull of $K$, is all of $\R^d$. Otherwise, we could restrict ourselves to the span of the $v_i$, on which $g$ still satisfies all properties of a norm.
In particular, we have $n \geq d$ and $\conv \{K \cup (-K)\} \in \K^d_o$, so $g$ can be assumed to be a norm on $\R^d$.
\end{remark}

\begin{theorem}
For any non-empty compact convex set $K \subset \R^d$ and any $u_1, \ldots, u_n \in \R^d$,
\[
    \max_{\varepsilon \in \{\pm1\}^n} h_K \Big( \sum_{i=1}^n \varepsilon_i u_i \Big)
    \geq \max \Big\{ \frac{1}{n}, \frac{1}{d} \Big\} \cdot \sum_{i=1}^n |h_K(u_i)|.
\]
The inequality is sharp for any pair $d,n \geq 1$.
\end{theorem}
\begin{proof}
To prove the inequality, we may assume that $h_K$ is a norm and $n \geq d$ by Remark~\ref{rem:sublin_reduction}.
The approach in Remark~\ref{rem:asymptotic_lower} with the first equality in \eqref{eq:1norm_ineq} omitted (since we do not necessarily have unit vectors) establishes the claimed inequality.
For the sharpness, take $K = B_1^d$ so that $h_K = \|.\|_\infty$ and let $\{e_j\}_{j=1}^d$ be the standard orthonormal basis of $\R^d$.
If $n \leq d$, choose $u_i = e_i$ for all $i \in [n]$.
Otherwise, choose $u_i = (n+1-d) e_i$ for $i \in [d-1]$ and $u_d = \dots = u_n = e_d$.
\end{proof}

\begin{theorem}
For any non-empty compact convex set $K \subset \R^d$ and any $u_1, \ldots, u_n \in \R^d$ with $|h_K(u_i)| \geq 1$ for all $i \in [n]$,
\[
    \max_{\varepsilon \in \{\pm1\}^n} h_K \Big( \sum_{i=1}^n \varepsilon_i u_i \Big)
    \geq \bigg\lceil \frac{n}{d} \bigg\rceil.
\]
The inequality is sharp for any pair $d,n \geq 1$.
\end{theorem}
\begin{proof}
Like before, we may assume that $h_K$ is a norm, where \eqref{eq:max_sublin} guarantees that all $u_i$ have norm at least $1$ after the reduction described in Remark~\ref{rem:sublin_reduction}.
The theorem now follows from the proof of Theorem~\ref{thm:sharplower} with the last equality in \eqref{eq:convexity_ineq} replaced by an inequality.
The sharpness follows from Theorem~\ref{thm:sharplower}.
\end{proof}

Our second (and final) generalization concerns the circumradius of Minkowski sums of compact convex sets. Lower estimates on such quantities 
were first studied in \cite{2012GonzalezHernandez} and have since been generalized to various other settings (cf.~\cite{2018AHS,2016CYL,2014Gonzalez,2014GonzalezHernandez,2024LombardiSaorin}), typically only for sums of two convex bodies.
See also \cite{2018BFGK}, where $n$-term sums
are considered in connection with vector sum problems.
Here, we establish a sharp lower bound on the Euclidean circumradius and diameter of Minkowski sums of compact convex sets with an arbitrary number of summands.
In fact, this result is an extension of Corollary~\ref{cor:eucl_general_vectors}, as is seen by taking $K_i = [-u_i,u_i]$ and referring to \eqref{eq:zonotope_LSS}.

\begin{theorem}
For any non-empty compact convex sets $K_1, \ldots, K_n \subset \R^d$,
\[
    R \Big( \sum_{i=1}^n K_i \Big)
    \geq \frac{2 \, \kappa_{d-1}}{d \, \kappa_d} \sum_{i=1}^n R(K_i),
\]
with equality if and only if $d=1$ or all $K_i$ are singletons. Moreover, the inequality is sharp for any $d \geq 1$ as $n \to \infty$, even under the extra assumption that all sets $K_i$ are symmetric.

All statements remain true for the diameter $D$ in place of the circumradius $R$.
\end{theorem}
\begin{proof}
We employ the argument of the proof of the lower bound in Proposition~\ref{prop:polball}: in \eqref{eq:zonotope_radius}, the sets $[-u_i,u_i]$ are replaced by the $K_i$, and $Q$ is set to be $\sum_{i=1}^n K_i$. The argument is still valid since $V_1$ is Minkowski linear on the space of non-empty compact convex sets. 

For the diameter, we simply note that for any non-empty compact convex set $K \subset \R^d$,
\[
    D(K)
    = R(K-K)
\]
and that our construction verifying the sharpness of the estimates only 
employs symmetric sets.
\end{proof}

\section*{Acknowledgements}
The authors are indebted to Imre Bárány, Tomasz Kobos, Szilárd Révész, Carsten Schütt and Gábor Tardos for the fruitful conversations and suggestions.

\bibliographystyle{siam}
\bibliography{Bibliography}

\bigskip

\noindent
{\sc Gergely Ambrus}
\smallskip

\noindent
{\em Bolyai Institute, University of Szeged, Hungary \\ and\\ 
 Alfréd Rényi Institute of Mathematics, Budapest, Hungary
 }
\smallskip

\noindent
e-mail address: \texttt{ambrus@server.math.u-szeged.hu}

\bigskip

\noindent
{\sc Florian Grundbacher}
\smallskip

\noindent
{\em Technical University of Munich}
\smallskip

\noindent
e-mail address: \texttt{florian.grundbacher@tum.de}
\end{document}